\documentclass{article}[12pt]

\usepackage{amsmath}
\usepackage{amssymb}
\usepackage{amsfonts}
\usepackage{amsthm}
\usepackage{setspace}
\usepackage{bbm}
\usepackage[english]{babel}
\usepackage{graphicx}
\usepackage{color}
\usepackage{fancyhdr} 
\usepackage[left=1in,right=1in,top=1in,bottom=1in]{geometry}
\usepackage{amsmath,amscd}
\usepackage[all]{xy}
\usepackage{tikz}
\usepackage{tikz-cd}
\usetikzlibrary{shapes}
\usetikzlibrary{arrows}
\usetikzlibrary{matrix}
\usetikzlibrary{calc}
\usetikzlibrary{positioning}
\usetikzlibrary{decorations.pathreplacing}
\usetikzlibrary{decorations.markings}
\usepackage{mathtools}
\usepackage[hypcap=true]{caption}
\usepackage{hyperref}
\usepackage[numbers]{natbib}
\DeclareFontFamily{U}{mathx}{\hyphenchar\font45}
\DeclareFontShape{U}{mathx}{m}{n}{
      <5> <6> <7> <8> <9> <10>
      <10.95> <12> <14.4> <17.28> <20.74> <24.88>
      mathx10
      }{}
\DeclareSymbolFont{mathx}{U}{mathx}{m}{n}
\DeclareFontSubstitution{U}{mathx}{m}{n}
\DeclareMathAccent{\widecheck}{0}{mathx}{"71}

\DeclareMathOperator{\Sp}{span}

\DeclareMathOperator{\GL}{GL}
\DeclareMathOperator{\SL}{SL}

\DeclareMathOperator{\Spec}{Spec}

\DeclareMathOperator{\diag}{diag}
\DeclareMathOperator{\ord}{ord}
\DeclareMathOperator{\cmid}{mid}
\DeclareMathOperator{\up}{up}
\DeclareMathOperator{\can}{can}
\DeclareMathOperator{\Conf}{Conf}

\DeclareMathOperator{\Cox}{Cox}
\DeclareMathOperator{\Pic}{Pic}
\DeclareMathOperator{\Fl}{\mathcal{F}\hspace{-1.6pt}\ell}
\DeclareMathOperator{\Fld}{\widetilde{\Fl}}

\newcommand{\tf}{\vartheta}
\newcommand{\mumgs}{\mu_{\mathrm{MGS}}}
\newcommand{\tc}[2]{{\textcolor{#1}{#2}}}
\newcommand{\lrp}[1]{\left(#1\right)}
\newcommand{\lrb}[1]{\left[#1\right]}
\newcommand{\lrm}[1]{\left|#1\right|}
\newcommand{\lrc}[1]{\left\{#1\right\}}
\newcommand{\lra}[1]{\langle{#1}\rangle}

\newcommand{\RA}{\Rightarrow}

\newcommand{\eq}[2]{\begin{equation}\label{#2} \begin{split} #1  \end{split} \end{equation}}
\newcommand{\eqn}[1]{\begin{equation*} \begin{split} #1 \end{split} \end{equation*}}

\newcommand{\R}{\mathbb{R} }
\newcommand{\C}{\mathbb{C} }

\newcommand{\Z}{\mathbb{Z} }

\newcommand{\kk}{\mathbbm{k} }
\newcommand{\gv}{\mathbf{g} }
\newcommand{\cv}{\mathbf{c} }
\newcommand{\vb}[1]{\mathbf{#1}}
\newcommand{\cA}{\mathcal{A} }
\newcommand{\cAp}{\mathcal{A}_{\mathrm{prin}} }

\newcommand{\prin}{{\mathrm{prin}} }
\newcommand{\cX}{\mathcal{X} }
\newcommand{\clV}{\mathcal{V} }
\newcommand{\cyU}{\mathcal{U} }
\newcommand{\ssO}{\mathcal{O} }

\newcommand{\cta}{\Conf_3(\Fld) }
\newcommand{\ctb}{\Conf_3(\Fl) }
\newcommand{\ctag}{\Conf_3^\times(\Fld) }
\newcommand{\ctbg}{\Conf_3^\times(\Fl) }
\newcommand{\ctal}{\widetilde{\Conf_3}(\Fld) }
\newcommand{\ctbl}{\widetilde{\Conf_3}(\Fl) }
\newcommand{\ctagl}{\widetilde{\Conf_3}^\times(\Fld) }
\newcommand{\ctbgl}{\widetilde{\Conf_3}^\times(\Fl) }
\newcommand{\Wl}{\widetilde{W} }
\newcommand{\Xil}{\widetilde{\Xi} }
\newcommand{\lb}{\mathcal{L} }
\newcommand{\trop}{\mathrm{trop} }
\newcommand{\uf}{\mathrm{uf} }

\newcommand{\abc}{\alpha, \beta, \gamma }
\newcommand{\Vabc}{V_\alpha \otimes V_\beta \otimes V_\gamma }
\newcommand{\itri}[1]{\vb{i}_{\triangle_{#1}} }
\newcommand{\trir}[1]{\triangle_{#1} }

\def\presuper#1#2%
  {\mathop{}%
   \mathopen{\vphantom{#2}}^{#1}%
   \kern-\scriptspace%
   #2}
\makeatletter   
\def\blfootnote{\xdef\@thefnmark{}\@footnotetext}
\makeatother
 
\theoremstyle{plain}
\newtheorem{theorem}{Theorem}
\newtheorem{prop}[theorem]{Proposition}
\newtheorem{lemma}[theorem]{Lemma}

\newtheorem{cor}[theorem]{Corollary}

\theoremstyle{definition}
\newtheorem{definition}[theorem]{Definition}

\theoremstyle{remark}
\newtheorem{remark}[theorem]{Remark}

\newenvironment{problem}[1]{ \flushleft \textcolor{blue}{\normalsize {#1}}}

\newenvironment{subprob}[1]{ \flushleft \textcolor{blue}{\normalsize {#1}}}

\mathtoolsset{centercolon}
\begin{document}

\renewcommand{\thesubsection}{\arabic{section}.\alph{subsection}}
\renewcommand{\labelenumi}{(\arabic{enumi})}

\title{Littlewood-Richardson coefficients via mirror symmetry for cluster varieties}
\author{Timothy Magee}
\begin{NoHyper}
\end{NoHyper}
\date{}
\maketitle

{\abstract{%
I prove that the full Fock-Goncharov conjecture holds for $\ctag$-- the configuration space of triples of decorated flags in generic position.
As a key ingredient of this proof, I exhibit a maximal green sequence for the quiver of the initial seed.
I compute the Landau-Ginzburg potential $W$ on $\ctag^\vee$ associated to the partial minimal model $\ctag \subset \cta$.
The integral points of the associated ``cone'' ${\Xi:=\lrc{W^T \geq 0} \subset \ctag^\vee\lrp{\R^T}}$
parametrize a basis for $\ssO\lrp{\cta} = \bigoplus \lrp{\Vabc}^G$
and encode the Littlewood-Richardson coefficients $c^\gamma_{\alpha \beta}$.
In the initial seed, the inequalities defining $\Xi$ are exactly the tail positivity conditions of \cite{Zel_Tail}.
I exhibit a unimodular $p^*$ map that identifies $W$ with the potential of Goncharov-Shen on $\ctag$\cite{GShen} and $\Xi$ with the Knutson-Tao hive cone.\cite{KTv1}
}}

\tableofcontents
\section{Introduction}
\subsection{\label{subsec:summary}Summary of results}
In this paper I obtain polytopes whose number of integral points are the Littlewood-Richardson coefficients by a method that has essentially nothing to do with representation theory.  
The same method will, in theory, produce analogous polytopes whose integral points parametrize a canonical basis for the space of sections of any line bundle on any Fano variety with a choice anti-canonical divisor.  
I prove Corollary 0.21 of \cite{GHKK}, and recover Corollary 0.20 as well with very little additional work.\footnote{These corollaries in \cite{GHKK} reference the current work and were stated as conjectures until a preliminary version of this paper was prepared.}
In particular I recover polytopes parametrizing canonical bases for each irreducible representation of $\GL_n$.

Let $G=\GL_n$, let $B$ be the Borel subgroup of upper triangular matrices in $G$, and let $U$ be the unipotent radical of $B$-- the subgroup of upper-triangular matrices with diagonal entries all $1$.
The flag variety $\Fl$ is isomorphic to $G/B$, and the {\it{decorated flag variety}} $\Fld$-- whose points consist of a complete flag ${X_\bullet= \lrp{X_1 \subset \cdots \subset X_n}}$ together with a non-zero vector $x_i$ in each successive quotient $X_i/X_{i-1}$-- is isomorphic to $G/U$.\footnote{We say these spaces are isomorphic rather than saying, {\it{e.g.}}, ``The flag variety {\emph{is}} $G/B$'' because no choice of Borel subgroup is necessary to describe the flag variety.  The flag variety is isomorphic to $G/B$ for any choice of Borel subgroup $B$, and we view $\Fl$ and $\Fld$ as being free of such choices.}
Setting $H:= B/U$, $\Fld$ is naturally a principal $H$-bundle over $\Fl$.
Next, following \cite{FG_Teich,GShen}, define $\ctb$ to be $ G \backslash \Fl^{\times 3}$
and $\ctbg$ to be the affine subvariety where pairs of flags intersect generically.
Define $\cta$ and $\ctag$ analogously-- $ \cta:= G \backslash \Fld^{\times 3}$ and
$\ctag \subset \cta$ is the generic locus.
Just like $\Fld \to \Fl$ is a principal $H$-bundle, 
$\cta$
is a principal $H^{\times 3}$-bundle over  
$\ctb$.
Moreover, in both cases, the base is Fano and the total space is the universal torsor for the base.
This $H^{\times 3}$-bundle has a very special property: 
\eqn{
\ssO\lrp{\cta} = \Cox\lrp{\ctb} := \bigoplus_{\lb \in \Pic\lrp{\ctb}}\Gamma\lrp{\ctb, \lb}.
}
Each irreducible representation of $G$ can be described as the space of sections of a line bundle over $\Fl$.
Similarly, for $\lb \in \Pic\lrp{\ctb}$, 
\eqn{\Gamma\lrp{\ctb, \lb} = \lrp{\Vabc}^G}
for some triple of weights $\lrp{\abc}$,
where $V_\lambda$ denotes the irreducible representation of highest weight $\lambda$.
The space $\cta$ is of interest here because
\eqn{\ssO\lrp{\cta} = \bigoplus_{\abc}\lrp{\Vabc}^G,}
and the dimension of each summand is a Littlewood-Richardson coefficient.
See Section~\ref{subsec:RT} for details.

The first key result of this paper pertains to the open subvariety $\ctag$.
{\theorem{\label{thm:FG}The full Fock-Goncharov conjecture holds for $\ctag$.}}

The {\it{full Fock-Goncharov conjecture}} \cite[Definition~0.6]{GHKK} will be described in detail in the next subsection~\ref{subsec:FG_intro}.
For now, there is a piecewise linear manifold $\ctag^\vee\lrp{\R^T}$ encoding the logarithmic geometry of the mirror $\ctag^\vee$ to $\ctag$.
Theorem~\ref{thm:FG} says that its integral points $\ctag^\vee\lrp{\Z^T}$ parametrize a canonical basis for $\ssO\lrp{\ctag}$.

Next observe that $\cta$ is a partial compactification of $\ctag$.
This compactification corresponds to a Landau-Ginzburg potential $W$ on $\ctag^\vee$.
A piecewise linear analogue of $W$-- its tropicalization $W^T$-- yields a subset $\Xi := \lrc{W^T \geq 0}$ of $\ctag^\vee\lrp{\R^T}$, and so a subset of our basis for $\ssO\lrp{\ctag}$.
Now, both $\ctag$ and $\ctag^\vee$ have an atlas of algebraic tori, with tori in the two atlases coming in dual pairs.
These tori are glued via birational maps, known as {\it{mutations}}.
Choosing a pair of tori identifies $\ctag^\vee\lrp{\R^T}$ with a real vector space and $\Xi$ with a rational polyhedral cone in this vector space.
Under one particularly simple choice of tori, the inequalities defining $\Xi$ are precisely the tail-positivity conditions of \cite{Zel_Tail}.\footnote{Thanks to G. Fourier and an anonymous referee for pointing this out to me.}

We need one more ingredient to identify $\Xi$ with the Knutson-Tao hive cone of \cite{KTv1}.
There is a class of maps (known as $p^*$ maps) that in this setting send the cocharacter lattice $N$ of a torus in $\ctag$ to the cocharacter lattice $M$ of the dual torus in $\ctag^\vee$ and commute with mutation.
Such a map of lattices gives rise to a map of schemes $p:\ctag \to \ctag^\vee$.\footnote{A complete discussion of $p^*$ maps appears in \cite[Section~2]{GHK_birational}, and I give a brief review in Section~\ref{subsec:WXip}.}
{\theorem{For a particular pair of dual tori,\footnote{For readers already familiar with the cluster structure of $\ctag$, these are the tori of the initial seed.} the cone $\Xi$ defined by the Landau-Ginzburg potential $W$ on $\ctag^\vee$ is precisely the tail-positive cone of \cite{Zel_Tail}.
Furthermore, a particular choice of the map $p^*:N \to M$ identifies $W$ with the potential of Goncharov-Shen on $\ctag$ and $\Xi$ with the Knutson-Tao hive cone.}}

Note that these identifications are not used to prove any of the combinatorial results in this paper. 
Instead, the identifications should be interpreted as placing the tail-positive cone and the Knutson-Tao hive cone within the broad framework of log Calabi-Yau mirror symmetry.
It is an instance of the mirror construction described, {\it{e.g.}}, in \cite{GHKK}.

\sloppy The mirror $\ctag^\vee$ to $\cta$ comes with a map to the dual torus $\lrp{H^{\times 3}}^\vee$,
tropicalizing to a map $$w: \ctag^\vee\lrp{\Z^T} \to \lrp{H^{\times 3}}^\vee\lrp{\Z^T}.$$
The integral tropicalization of a torus $T$ is just its cocharacter lattice $\chi_*\lrp{T}$, so
$\lrp{H^{\times 3}}^\vee\lrp{\Z^T}$ is the character lattice $\chi^*\lrp{H^{\times 3}}$.
Theorem~\ref{thm:FG}, 
together with the existence of an {\it{optimized seed for each frozen variable}}\footnote{This condition allows us to associate a regular function on $\ctag^\vee$ to each irreducible component of ${D:=\cta \setminus \ctag}$, and it ensures $D$ interacts well with the canonical basis for $\ssO\lrp{\ctag}$. This will be discussed in greater detail in Section~\ref{subsec:optimized}.} (Proposition~\ref{prop:op}) 
and existence of a unimodular $p^*$ map (Proposition~\ref{prop:unimodular}), 
implies that points in ${\Xi\lrp{\Z^T}:= \Xi \bigcap \ctag^\vee\lrp{\Z^T}}$ are canonically identified with regular functions on $\cta$ invariant under the $H^{\times 3}$ action.
If $q \in  \Xi \lrp{\Z^T}$ and $\tf_q$ is the corresponding function on $\cta$, then $w\lrp{\tf_q}$ is the weight of $\tf_q$ under the $H^{\times 3}$ action.
Given a weight $\lrp{\alpha,\beta,\gamma}$ of this action, $w^{-1}\lrp{\alpha,\beta,\gamma} \bigcap \Xi\lrp{\Z^T}$ parametrizes a basis for the $\lrp{\alpha,\beta,\gamma}$-weight space of $\ssO\lrp{\cta}$.
Since $\ssO\lrp{\cta} = \bigoplus_{\alpha,\beta,\gamma} \lrp{V_\alpha \otimes V_\beta \otimes V_\gamma}^G$, counting these points gives the Littlewood-Richardson coefficients.  
This is described in more detail in Section~\ref{sec:Conf3back}.

{\remark{$\cta$ and $\Fld$ have very similar cluster structures.  As a result, many of the proofs in \cite{Magee} apply here as well.  For completeness and convenience, I have provided them here.  
Since the current paper encompasses the main results in \cite{Magee}, I'll seek publication of this paper and not \cite{Magee}.}}

\subsection{\label{subsec:FG_intro}Full Fock-Goncharov conjecture}

This subsection provides a bit of background on the full Fock-Goncharov conjecture following \cite{GHKK}.

The full Fock-Goncharov conjecture is a canonical basis conjecture for regular functions on {\it{cluster varieties}}.
Cluster varieties are schemes built out of birationally glued algebraic tori.
A torus $\lrp{\C^*}^n$ comes with a canonical volume form $\frac{dz_1}{z_1}\wedge \cdots \wedge \frac{dz_1}{z_1}$.
The birational gluing maps are required to patch these volume forms on tori, giving a global volume form.
In cluster theory, there are two classes of such gluing maps ({\it{mutations}}), producing two types of cluster varieties-- known as $\cA$-varieties and $\cX$-varieties.
The $\cA$ and $\cX$ mutations are defined dually from the same data, and up to certain multipliers $\cA$ and $\cX$ are built from dual tori.
There is a notion of {\it{Langlands dual data}} that accounts for these multipliers, giving rise to another pair of cluster varieties $\lrp{\presuper{L}{\cA},\presuper{L}{\cX}}$ built precisely out of the dual tori for $\lrp{\cA,\cX}$.
The pairs $\lrp{\cA,\presuper{L}{\cX}}$ and similarly $\lrp{\presuper{L}{\cA},\cX}$ are said to be {\it{Fock-Goncharov dual}}.
For a discussion of $\cA$ and $\cX$ varieties and their Fock-Goncharov duals, please see \cite[Section~1.2]{FG_cluster_ensembles}, \cite[Section~2]{GHK_birational}, and \cite[Appendix~A]{GHKK}.


Let $\clV$ be a cluster variety (of whichever type), and $\clV^\vee$ its Fock-Goncharov dual.
We think of $\clV^\vee$ as mirror to $\clV$ as, if $\clV$ and $\clV^\vee$ are sufficiently close to affine, the Fock-Goncharov canonical basis conjecture becomes a special case of the log Calabi-Yau mirror symmetry conjecture \cite[Conjecture~0.6]{GHK_logCY}.
This is discussed in greater detail in the introduction of \cite{GHKK},
where they additionally give a log Calabi-Yau mirror symmetry interpretation of the Fock-Goncharov canonical basis conjecture for arbitrary cluster varieties-- not only those that satisfy certain affineness assumptions.\footnote{While the original Fock-Goncharov conjecture is false as stated (\cite{GHK_birational}), \cite{GHKK} build toric degenerations of $\cA$-varieties and show that the conjecture holds in a formal neighborhood of the central fiber of the degeneration.  
They interpret this as saying the Fock-Goncharov conjecture always holds in the {\it{large complex structure limit}}.}  

In \cite{GHKK}, several algebras are associated to $\clV$.
First, there is the {\it{upper cluster algebra}}
$\up(\clV) = \ssO\lrp{\clV}$, originally defined for $\cA$-type cluster varieties in \cite{FZ_clustersIII}.
Its subalgebra generated by global monomials,
{\it{i.e.}} global regular functions restricting to a character on some torus in the atlas for $\clV$,
is the {\it{ordinary cluster algebra}} $\ord(\clV)$.
In the case of an $\cA$-type cluster variety, this corresponds to the usual notion of a cluster algebra.
If $\clV^\vee$ is the Fock-Goncharov dual of $\clV$,
$\can(\clV)$ is a vector space with basis parametrized by $\clV^\vee\lrp{\Z^T}$.\footnote{This set-- the {\it{integral tropical points of $\clV^\vee$}}-- is essentially the divisors $D$ that we can use to partially compactify $\clV^\vee$ such that the canonical volume form on $\clV^\vee$ has a pole along $D$.  It will be discussed further in the following subsection~\ref{subsec:Compact}. For a more complete description of tropicalization, see \cite[Section~2]{GHKK}.}
{\it{Scattering diagrams}} and {\it{broken lines}} are used to associate
to each $q \in \clV^\vee\lrp{\Z^T}$ a (possibly infinite) sum of characters on each torus in $\clV$'s atlas,
the result denoted by $\vartheta_q$,
and to define a multiplication rule for the $\vartheta_q$.
A scattering diagram for $\clV$ is a collection of walls, each decorated with a scattering function, living in a piecewise linear manifold $\clV^\vee\lrp{\R^T}$.
Broken lines are lines drawn in this piecewise linear manifold which are allowed to bend in a prescribed fashion at the walls.
Each linear piece is decorated with a monomial corresponding to the tangent direction, and the allowed bending at a wall is determined by the wall's scattering function.
The $\vartheta$-functions are expressed in patches by summing decorating monomials over broken lines with a given initial direction and end point. 
The details of this construction are beyond the scope of this paper.
In situations where the full Fock-Goncharov conjecture holds,
$\can(\clV)$ will be identified with $\up(\clV)$,
and the $\vartheta_q$ will form a canonical basis for $\up(\clV)$.
More generally, $\can(\clV)$ has a subspace $\cmid(\clV)$ parametrized by the subset $\Theta \subset \clV^\vee\lrp{\Z^T}$
consisting of $\vartheta_q$ which
restrict to finite sums of characters, {\it{i.e.}} Laurent polynomials, on tori from the atlas.
Then each element of $\cmid(\clV)$ naturally corresponds to an element of $\up(\clV)$,
but there is no reason {\it{a priori}} that distinct elements of $\cmid(\clV)$ must correspond to distinct elements of $\up(\clV)$.
More formally, there is a canonical algebra homomorphism $\nu: \cmid(\clV) \to \up(\clV)$, and we have the following definition.

{\definition{\cite[Definition~0.6]{GHKK} We say that the {\it{full Fock-Goncharov conjecture holds for $\clV$}} if $\nu: \cmid(\clV) \to \up(\clV)$ is injective, $\up(\clV) = \can(\clV)$, and $\Theta = \clV^\vee\lrp{\Z^T}$.}}

To recap, we have the following: 
\begin{enumerate}
    \item \label{up} $\up(\clV) = \ssO\lrp{\clV}$
    \item \label{ord} $\ord(\clV) \subset \up(\clV)$, the subalgebra generated by global monomials
    \item \label{can} $\can(\clV)$, the vector space spanned by $\vartheta_q$ for $q\in \clV^\vee\lrp{\Z^T}$, endowed with a multiplication rule via scattering diagrams and broken lines
    \item  \label{mid} $\cmid(\clV) \subset \can(\clV)$, the span of $\Theta\subset \clV^\vee\lrp{\Z^T} $, the $\vartheta_q$ which are Laurent polynomials on cluster tori.
\end{enumerate}
The full Fock-Goncharov conjecture pertains to (\ref{up}), (\ref{can}), and (\ref{mid}), rather than the usual notion of a cluster algebra (\ref{ord}). The statement that the full Fock-Goncharov conjecture holds for $\clV$ is essentially that (\ref{up}), (\ref{can}), and (\ref{mid}) are all equal.
 
Many conditions implying the full Fock-Goncharov conjecture holds for a given cluster variety are provided in \cite{GHKK}.
I will use 

{\bf{\cite[Proposition~8.25]{GHKK}:}}  {\it{If $\cA$ has large cluster complex, then $\cAp$ has Enough Global Monomials, $\Theta = \cAp^\vee\lrp{\Z^T}$, and the full Fock-Goncharov conjecture holds for $\cAp$,
$\cX$, very general $\cA_t$, and, if the convexity condition (7) of \cite[Theorem~0.3]{GHKK} holds, for $\cA$. }}

The {\it{cluster complex}} is a subset of the scattering diagram for $\cA$, where chambers correspond to tori in the atlas for $\cA$ and rays spanning these chambers correspond to cluster variables.\footnote{See \cite[Definition~2.9]{GHKK} for the precise definition of the cluster complex.}
Broken lines initiating within this region correspond to {\it{cluster monomials}}-- $\tf$-functions that restrict to characters on some torus of $\cA$'s atlas.
Each choice of torus in $\cA$'s atlas gives an identification of $\cA^\vee\lrp{\R^T}$ with a real vector space.  
The statement that $\cA$ has {\it{large cluster complex}} means that, for some choice of torus, the cluster complex is not contained in a half-space.

Before summarizing how I prove that $\ctag$ has large cluster complex, I need to provide a bit more background.
Cluster varieties are built recursively from a $\Z$-basis for the character lattice of a torus and data that indicates how to replace this torus and basis with new tori and bases.
The bases are known as {\it{seeds}}, and the data used to replace one seed with another is often encoded in a {\it{quiver}}-- a directed graph.
Mutation changes both the seed and the quiver.

To show that $\ctag$ has large cluster complex, I exhibit a {\it{maximal green sequence}} \cite[Definition~2.8]{BDP} for the quiver of the initial seed.
I review the notion of a maximal green sequence in Section~\ref{subsec:MGS}.
The existence of a maximal green sequence implies $\ctag$ has large cluster complex by \cite[Proposition~8.24]{GHKK}.\footnote{The weaker condition of a {\it{reddening sequence}}-- see {\it{e.g.}} \cite{Keller_DT}-- would also suffice here.}

{\theorem{\label{thm:green}The quiver for the initial seed of $\ctag$ has a maximal green sequence, and therefore $\ctag$ has large cluster complex.}}

Next, the convexity condition (7) of \cite[Theorem~0.3]{GHKK} referenced in \cite[Proposition~8.25]{GHKK} is the following:

\quad {\it{There is a seed $\vb{s}= \lrp{e_1, \dots, e_n}$ for which all the covectors $\lrc{e_i, \cdot}$, $i\in I_{\uf}$, lie in a strictly convex cone.}}  

Here, $\vb{s}$ is a basis for the cocharacter lattice $N$ of a torus in the $\cA$-variety.
The data used to define $\cA$ includes a skew form $\lrc{ \cdot, \cdot}$ on $N$, an indexing set $I$ of cardinality $n$, and a subset $I_{\uf}\subset I$ corresponding to allowed directions of mutation. 
In writing {\it{covectors $\lrc{e_i, \cdot}$, $i\in I_{\uf}$}}, Gross-Hacking-Keel-Kontsevich refer to this skew form and subset of $I$. 

I show that the convexity condition holds for the initial seed.
Together with Theorem~\ref{thm:green}, this shows Theorem~\ref{thm:FG}-- the full Fock-Goncharov conjecture holds for $\ctag$.

\subsection{\label{subsec:Compact}Partial compactifications and potentials}

The space we are really interested in is $\cta$, rather than $\ctag$.
It is $\cta$, not $\ctag$, that gives the decomposition
\eqn{
\ssO\lrp{\cta} = \bigoplus_{\alpha, \beta, \gamma} \lrp{V_\alpha \otimes V_\beta \otimes V_\gamma}^G.
}
But Subsection~\ref{subsec:FG_intro} describes a canonical basis for $\ssO\lrp{\ctag}$.

This situation is typical.
Generally spaces we are interested in, say for representation theoretic reasons, will not be cluster varieties, or even log Calabi-Yau varieties.\footnote{%
Cluster varieties, both $\cA$- and $\cX$-varieties, form a very special class of log Calabi-Yau varieties, 
which generalize Calabi-Yau varieties to the non-compact world. 
See \cite[Definition~1.1]{GHK_birational}.  
I frame this discussion in terms of log Calabi-Yau varieties rather than just cluster varieties in this introductory section simply because it is the log Calabi-Yau structure, not the cluster structure, that is essential to the narrative.
The cluster structure will be very useful in many of the proofs later on, but this is a case where framing the narrative in a more general context clarifies the key arguments rather than obscuring them.}
However, many representation theoretically interesting spaces {\it{are}} partial compactifications of cluster varieties or log Calabi-Yau varieties in a nice way.
For log Calabi-Yau varieties, the ``nice'' type of partial compactification we're interested in is called a {\it{partial minimal model}}.
Take $\cyU$ to be a log Calabi-Yau with canonical volume form $\Omega$.
Then an inclusion $\cyU \subset Y$ as an open subset is a {\it{partial minimal model}} if $\Omega$ has a simple pole along every irreducible divisor of $Y\setminus \cyU$.
In the special case that $\cyU$ is a cluster $\cA$-variety with frozen variables, there is a simple way these partial minimal models may arise-- 
by taking $Y$ to be the partial compactification given by allowing some frozen variables to vanish.\cite[Section~0.3]{GHKK}
$\ctag$ and $\cta$ are related in precisely this way.

In this situation, each irreducible divisor in $D:= Y\setminus \cyU$ defines a {\it{divisorial discrete valuation}} whose evaluation at $\Omega$ is negative.
So, each irreducible component gives a point in $\cyU^\trop\lrp{\Z}$ {\emph{by definition}}.\cite[Definition~1.7]{GHK_birational}
A {\it{divisorial discrete valuation}} (ddv) is a discrete valuation $v: \C\lrp{\cyU} \setminus \lrc{0} \to \Z$ on the field of rational functions of $\cyU$ given by order of vanishing along a divisor on some variety birational to $\cyU$.
The set $\cyU^\trop\lrp{\Z}$ is defined to be 
\eqn{\lrc{\left.v: \C\lrp{\cyU} \setminus \lrc{0} \to \Z\right| v \text{ is a ddv, } v\lrp{\Omega} < 0  } \cup\lrc{0}.}
This definition applies to all log Calabi-Yau varieties.

There is a related notion that specifically references the atlas of tori and mutation maps, known as {\it{Fock-Goncharov tropicalization}}.
Fock-Goncharov tropicalization applies to {\it{positive schemes}}.
Let $\cyU$ be a scheme where an open dense subset has an atlas of tori, glued birationally.
Suppose additionally that the transition functions between tori have the form $\frac{f}{g}$ where $f$ and $g$ are linear combinations of characters (Laurent monomials) with coefficients in $\Z_{\geq 0}$.
Such a $\cyU$ is said to be a {\it{positive scheme}}.\cite[Definition~1.1]{FG_cluster_ensembles}
If $P$ is a semifield, we can take the $P$-points of each torus-- $T\lrp{P}:= \chi_*\lrp{T} \otimes_\Z P$, where $\chi_*\lrp{T}$ is the cocharacter lattice of $T$ and we use the abelian group structure of $P$ for the tensor product.
The transition functions of a positive scheme induce isomorphisms of the $P$-points of the tori, giving a notion of the $P$-points of $\cyU$.
In particular, $P$ can be taken to be a {\it{tropical semifield}} such as $\Z^T$-- the set $\Z$ with multiplication $a\otimes b := a + b$ and addition $a \oplus b:= \max\lrp{a,b}$. 
This gives the notion of the Fock-Goncharov tropicalization $\cyU\lrp{\Z^T}$.

Cluster varieties are both log Calabi-Yau varieties and positive schemes, so both notions of tropicalization make sense here.
They are closely related as well.
There is a different tropical semifield structure on $\Z$, denoted $\Z^t$, given by replacing $\max$ with $\min$.
For a cluster variety $\cyU$, the notions $\cyU^\trop\lrp{\Z}$ and $\cyU\lrp{\Z^t}$ coincide.
We have an isomorphism of the tropical semifields $\Z^T$ and $\Z^t$ given by sign change $x \mapsto -x$. 
So, sign change also induces an identification $i: \cyU\lrp{\Z^T} \to \cyU^\trop\lrp{\Z}$.
These different notions of tropicalization and the relations between them are carefully discussed in \cite[Section~2]{GHKK}. 
With regard to the current paper, the reader is encouraged to ignore the distinction between the two forms of tropicalization, only paying it mind in the event that some sign confusion arises.
The identification map $i$ can be thought of as a minor detail that simply allows us to relate differing conventions in the literature.


From now on, take $\cyU$ to be a cluster variety.
Evaluation gives a pairing between $\cyU^\trop\lrp{\Z}$ and $\C\lrp{\cyU}$.
Since each element of $\Theta$ is canonically identified with a regular function on $\cyU$, we can restrict this pairing from $\C\lrp{\cyU}$ to $\Theta$.
We can similarly restrict from $\cyU^\trop\lrp{\Z}$ to $\Theta^\vee$.
In this way, evaluation defines a pairing
\eqn{
\lra{ \cdot, \cdot }: \Theta^\vee \times \Theta &\to \Z\\
			\lrp{ v, w } &\mapsto v\lrp{\vartheta_w}. 
}
We could just as well start with the evaluation pairing between $\lrp{\cyU^\vee}^\trop\lrp{\Z}$ and $\C\lrp{\cyU^\vee}$, which would restrict to the pairing
\eqn{
\lra{ \cdot, \cdot }^\vee: \Theta^\vee \times \Theta &\to \Z\\
			\lrp{ v, w } &\mapsto w\lrp{\vartheta_v}. 
}
%
These two pairings are conjectured to agree in general, and they are known to agree when either $v$ or $w$ corresponds to a global monomial.\cite[Lemma~9.10 and Remark~9.11]{GHKK}

Since the irreducible components $D_i$ of the divisor $D$ in our partial minimal model $Y \supset \cyU$ define points $v_i$ in $\cyU^\trop\lrp{\Z}$,
it's natural to ask whether these $v_i$ further lie in $\Theta^\vee$, canonically giving regular functions $\vartheta_{v_i}$ on $\cyU^\vee$.
If so, and if $\lra{ v_i , \cdot} = \lra{ v_i , \cdot}^\vee$, then
\eqn{\vartheta_{v_i}^\trop(w):= w\lrp{\vartheta_{v_i}} = v_i \lrp{\vartheta_w} = \ord_{D_i}\lrp{\vartheta_w} .}
So in this scenario,  $\vartheta_w$ extends to $D_i$ if and only if $\vartheta_{v_i}^\trop(w)\geq 0$, and it extends to $D$ if and only if
\eqn{ \min_i \lrc{ \ord_{D_i}\lrp{\vartheta_w} } = w\lrp{W} =: W^\trop(w) \geq 0 ,}
where
$W:=\sum_i \vartheta_{v_i}$.  
This function $W$ is called the {\it{Landau-Ginzburg potential}} associated to the partial minimal model $Y\supset \cyU$.

Assuming the full Fock-Goncharov conjecture holds for $\cyU$, this gives a candidate basis for $\ssO(Y)$, namely $\lrc{\vartheta_w \left| W^\trop \lrp{w} \geq 0\right.}$.
This need not be a basis though.
Poles can cancel when we add functions, so in principal we could have $\vartheta_p+\vartheta_q$ regular on $Y$ even if $\vartheta_p$ and $\vartheta_q$ have poles along $D$.

So, once we know that $\lrp{\cyU^\vee}^\trop\lrp{\Z}$ is canonically identified with a basis for $\ssO(\cyU)$, there are three issues we need to address in order to cut this down to a basis for $\ssO(Y)$.
\begin{enumerate}
    \item Each irreducible component $D_i$ of $D$ canonically determines a regular function $\theta_{v_i}$ on $\cyU^\vee$.
    \item For each component, $\lra{ v_i , \cdot} = \lra{ v_i , \cdot}^\vee$.
    \item If a linear combination of $\vartheta$-functions is regular on $Y$, then so is each $\vartheta$-function appearing in the sum.
\end{enumerate}

If $\cyU$ is a cluster $\cA$-variety and $Y$ is obtained by allowing frozen variables to vanish, all three hold if each frozen index has an optimized seed.
By construction, if a seed $s$ is optimized for the frozen index $i$, then $\vartheta_{v_i}$ restricts to a Laurent monomial on the cluster torus in $\cX$ associated to $s$.
I provide a brief review of optimized seeds in Section~\ref{subsec:optimized}. 
A more complete discussion can be found in \cite[Section~9.1]{GHKK}. 

In Section~\ref{subsec:optimized}, I show that each frozen index for $\ctag$ has an optimized seed.


%
{\theorem{Let $\ctag^\vee$ denote the Fock-Goncharov dual of $\ctag$, and let $W$ be the Landau-Ginzburg potential on $\ctag^\vee$ associated to the partial minimal model
$\ctag \subset \cta$.  Then
$$\Xi \lrp{\Z^T} := \lrc{W^T \geq 0 } \bigcap \ctag^\vee\lrp{\Z^T}$$
is canonically identified with a basis for $\ssO\lrp{\cta}$, determined by the pair $\ctag \subset \cta$.}}
 
In general if we have some partial minimal model $\cyU \subset Y$ for a log Calabi-Yau, we can't expect to get a canonical basis for $\ssO(Y)$ itself.
The basis will be determined by the geometry of the pair $\cyU \subset Y$ rather than $Y$'s geometry alone.
However, in the particular case $Y=\cta$, the log Calabi-Yau open subset $\cyU = \ctag$ is determined entirely by $Y$'s own geometry.
It is simply the locus where underlying flags intersect generically, described in more detail in Section~\ref{sec:Conf3back}.
In this sense, $\Xi\lrp{\Z^T}$ can be viewed as a canonical basis for $\ssO\lrp{\cta}$ itself-- it is a basis determined entirely by $\cta$'s own geometry.

{\it{Acknowledgements: }} 
First I would like to thank my advisor, Sean Keel.  He introduced me to the topics discussed here and proposed the problems I address.  I appreciate his guidance and poker games, and M. Moore's renowned soup.
Conversations with I. Ganev, S. Gunningham, and L. Shen cleared up a variety of confusions that arose in the course of this work.
I thank L. Bossinger for carefully reading through a preliminary version of this paper, and providing numerous edits and suggestions. 
I would also like to thank the two anonymous referees.  Both provided amazingly thorough and knowledgeable reports, and I feel their suggestions have greatly improved the exposition here.
During the preparation of this paper, I was partially supported by the NSF RTG grant DMS 1148490, by FORDECYT 265667, and by the EPSRC grant EP/P021913/1.

\section{\label{sec:Conf3back}Discussion of \texorpdfstring{$\cta$}{configurations of flags}}
\subsection{\label{subsec:RT}Representation theory background}

Interest in $\cta$ has its roots in representation theory.
The starting point is the Peter-Weyl theorem.
A group $G$ acts on itself both by left and right multiplication, and this action gives $\ssO\lrp{G}$ the structure of a $G\times G$-bimodule.
The following statement of the Peter-Weyl theorem comes from \cite{Procesi}.
\begin{theorem}
(Peter-Weyl) Let $G$ be a linearly reductive group.
Then as $G\times G$-bimodules
\eqn{
\ssO\lrp{G} = \bigoplus_{\lambda} V_\lambda \otimes V_\lambda^*,
}
where the sum is over isomorphism classes of irreducible representations of $G$.
\end{theorem}
For $\GL_n$, the span of the highest weight vector $v_\lambda$ in the irreducible representation $V_\lambda$ of highest weight $\lambda$
is the one dimensional subspace fixed by $U$-- the subgroup of upper triangular matrices with $1$'s along the diagonal.
So 
\eq{
\ssO\lrp{G}^{1\times U} = \bigoplus_{\lambda} V_\lambda \otimes \C\cdot u_\lambda,
}{eq:PW}
where $u_\lambda$ is the highest weight vector for $V_\lambda^*$.
The weight of $u_\lambda$ is $-w_0\lrp{\lambda}$, where $w_0$ is the longest element of the Weyl group $W$ of $G$.\footnote{The Weyl group for $\GL_n$ is just $S_n$, 
and $w_0$ is the permutation sending $1,2,\dots,n$ to $n,n-1,\dots,1$.}
So the copy of $V_\lambda$ appearing in Equation~(\ref{eq:PW}) is a weight space for the right action of the maximal torus $H$ in $G$, 
and its weight is $-w_0\lrp{\lambda}$.
To stress this point, the left action of $H$ splits $V_\lambda$ into weight spaces, the highest weight being $\lambda$, but under the right action $V_\lambda$ {\emph{is}} the $-w_0\lrp{\lambda}$-weight space. 
The next thing to observe is that functions on $G$ that are fixed by $U$-- so $f\lrp{x u} = f(x) $ for all $u \in U$-- are the same as functions on $G/U \cong \Fld$.
Then
\eqn{
\ssO\lrp{\Fld} = \bigoplus_\lambda V_\lambda,
}
and this is a weight space decomposition for the right action of $H$.

Now if we were to take three copies of $\Fld$ instead of one, we would have
\eqn{
	\ssO\lrp{\Fld^{\times 3}} = \bigoplus_{\alpha, \beta, \gamma} V_\alpha \otimes V_\beta \otimes V_\gamma.
}
The $G$-fixed subspace of $\Vabc$ (with $G$ acting on the left diagonally) is identified with $\hom_G\lrp{V_\gamma^*,V_\alpha \otimes V_\beta}$.
By Schur's lemma, this is just a copy of the trivial representation for every copy of $V_\gamma^*$ appearing in $V_\alpha \otimes V_\beta$.
Then
\eqn{
\ssO\lrp{\cta} = \bigoplus_{\abc}\lrp{\Vabc}^G,
}
with $\lrp{\Vabc}^G$ the $\lrp{-w_0\lrp{\alpha}, -w_0 \lrp{\beta}, -w_0 \lrp{\gamma}}$-weight space of the right $H^{\times 3}$ action,
and
\eqn{
\dim\lrp{\Vabc}^G = c^{-w_0\lrp{\gamma}}_{\alpha \beta}.
}
The term on the right is a {\it{Littlewood-Richardson coefficient}}.  These are the structure constants giving the decomposition 
$$V_\alpha \otimes V_\beta =  \bigoplus_\gamma V_\gamma^{c^\gamma_{\alpha \beta}}.$$ 
This is the connection between $\cta$ and the Littlewood-Richardson coefficients.

\subsection{\label{subsec:logCY}Geometric background}
As mentioned in Section~\ref{subsec:summary},
$\ctb$ is Fano and $\pi: \cta \to \ctb$ is naturally an $H^{\times 3}$-bundle.
Points in $\ctb$ are triples of complete flags, defined up to an overall $G$-action.
Given two arbitrary flags 
$X_\bullet = \lrp{X_1 \subset \cdots \subset X_n}$ 
and  
$Y_\bullet = \lrp{Y_1 \subset \cdots \subset Y_n}$ ,
we expect the $i$-dimensional subspace $X_i$ and the $\lrp{n-i}$-dimensional subspace $Y_{n-i}$ 
to intersect transversely.
A triple of flags $\lrp{X_\bullet,Y_\bullet,Z_\bullet} \in \ctb$ is in {\it{generic configuration}}
if each pairwise transversality condition is satisfied.\footnote{This is a restatement of the definitions given in \cite{FG_Teich,GShen}, specific to $\ctb$.}
$\ctbg \subset \ctb$ is the subset consisting of such triples of flags.
It is log Calabi-Yau--
its complement is an anticanonical divisor $D$ in $\ctb$. 
Furthermore, the canonical volume form on $\ctbg$ has a pole along all of $D$.
We could in principle use the log Calabi-Yau mirror symmetry machinery to study the pair $\ctbg \subset \ctb$.
After all, the vector spaces of interest $\lrp{\Vabc}^G$ from Subsection~\ref{subsec:RT} are spaces of sections of line bundles over $\ctb$.
However, lifting to $\cta$ will allow us to tackle all of the line bundles, and so all of the vector spaces $\lrp{\Vabc}^G$, at once.
$\ctag$ is precisely $\pi^{-1}\lrp{\ctbg}$ and $\cta$ is again a (partial) minimal model for $\ctag$.
 
%
%

\subsection{\label{subsec:cluster}Cluster structure}
$\ctag$ is not just log Calabi-Yau.
Fock and Goncharov described a cluster structure for it in \cite{FG_Teich}.
The discussion here is based on \cite{FG_Teich} and \cite{GShen}.

Define $\ctal := \SL_n \backslash \lrp{\GL_n/U}^{\times 3}$, and define $\ctagl$, $\ctbl$, and $\ctbgl$ analogously.
It will also be handy later to define $\Wl$ and $\Xil$ to be the Landau-Ginzburg potential on $\ctagl^\vee$ and the cone given by its tropicalization.
I'll describe the initial seed of $\ctagl$, viewed as a cluster $\cA$-variety, and we'll view $\ctag$ as a quotient of $\ctagl$.

The quiver for the initial seed comes from the ``$n$-triangulation'' of a triangle, illustrated below for $n=4$.

\input{4triangulation.tex}

The vertices in the $n$-triangulation will be the vertices of our quiver.
For the arrows, we need to orient the edges of the $n$-triangulation.
First, the boundary of the original triangle is given a clockwise orientation.  
The edges of the $n$-triangulation inherit their orientation from this one in the manner illustrated below.

\input{4triangulationOriented.tex}

The vertices on the boundary of the original triangle are frozen vertices of the quiver, and the vertices in the interior are unfrozen.
We ignore arrows between frozen vertices as they do not affect mutation, so the quiver we are after is

\input{InitialQuiver.tex}

The vertices of the quiver can be indexed by triples of non-negative integers $\lrp{a,b,c}$ satisfying $a+b+c = n$.

\input{IndexedQuiv.tex}

Take $V$ to be an $n$-dimensional vector space.
A point $X$ in $\Fl \cong \GL\lrp{V}/U\lrp{V}$ is a complete flag $X_\bullet = \lrp{X_1 \subset \cdots \subset X_n}$ of subspaces of $V$ 
together with non-zero vector $x_i$ in each successive quotient $X_i / X_{i-1}$.
I'll denote this by $x_\bullet = \lrp{x_1,\dots,x_n}$.
Now choose a volume form $\omega$ on $V$.
The cluster variables in the initial seed of $\ctagl$ are defined as follows:
\eqn{
A_{\lrp{a,b,c}}: \lrp{X,Y,Z} \mapsto \omega\lrp{x_1 \wedge \cdots \wedge x_a \wedge y_1 \wedge \cdots \wedge y_b \wedge z_1 \wedge \cdots \wedge z_c }.
}
Note that by definition a linear transformation
$T: V \to V$ 
is in $\SL\lrp{V}$ 
if and only if
$\wedge^n T$ acts by the identity on $\wedge^n V$.
So $A_{\lrp{a,b,c}}$ is indeed a well defined function on $\ctagl$-- it respects the quotient by the diagonal $\SL\lrp{V}$ action.
None of these cluster variables are invariant under the diagonal action of $\GL\lrp{V}$-- they are {\emph{not}} functions on $\ctag$--
but rational functions in these variables can still be $\GL\lrp{V}$ invariant.
In fact, take a Laurent monomial in these variables: $$f= \prod_{a+b+c = n} A_{\lrp{a,b,c}}^{r_{\lrp{a,b,c}}}.$$
Then for $g \in \GL\lrp{V}$, 
\eqn{
g\cdot f \lrp{X,Y,Z} = \lrp{\det g}^{\sum r_{\lrp{a,b,c}}} f \lrp{ X,Y,Z }.
}
So $f$ is $\GL\lrp{V}$ invariant if and only if
\eqn{
\sum_{a+b+c = n} r_{\lrp{a,b,c}} = 0.
}
This will lead to a condition on {\it{$\gv$-vectors}} (Proposition~\ref{prop:0sum}), which will index our $\vartheta$-functions.  I review the notion of $\gv$-vectors in Section~\ref{subsec:WXip}.
 
{\remark{\label{rem:CS}The initial data I have described is for $\ctagl$ rather than $\ctag$.  
That said, it can easily be translated into initial data for $\ctag$.  
That is, we can view $\ctag$ as a cluster variety in its own right, rather than as a quotient of a cluster variety.  
One way to do this is to replace all cluster variables $A_{\lrp{a,b,c}}$ of the initial seed with a new collection of variables, say 
$\overline{A}_{\lrp{a,b,c}} = A_{\lrp{a,b,c}}/ A_{\lrp{n,0,0}}$.
Upon doing so, the proofs I give in the following sections using $\ctagl$'s cluster structure translate immediately to $\ctag$.
However, I find it more natural to avoid such choices.
In what follows, I will freely use $\ctagl$'s cluster structure without further comment.}}

\section{\label{sec:FG}Full Fock-Goncharov conjecture holds for \texorpdfstring{$\ctag$}{generic configurations}}
I will show that the full Fock-Goncharov conjecture holds for $\ctag$ by proving the following two conditions:
\begin{enumerate}
	\item {\label{item:MGS}}The quiver $Q_{s_0}$ for the initial seed of $\ctag$ has a maximal green sequence.
	\item {\label{item:convex}}In the initial seed $s_0 = \lrp{e_1,\dots,e_n}$, all of the covectors $\lrc{ e_i, \cdot }$, $i\in I_{\mathrm{uf}}$, lie in a strictly convex cone.\footnote{$I_{\uf}$ is the unfrozen subset of the indexing set $I$-- the subset corresponding to mutable vertices.}
\end{enumerate}
Together, (\ref{item:MGS}) and (\ref{item:convex}) imply that the full Fock-Goncharov conjecture holds for $\ctag$.\cite[Proposition~8.25]{GHKK}
We'll begin with (\ref{item:convex}) as its proof is much shorter.
{\prop{\label{prop:convex}In the initial seed $s_0 = \lrp{e_1,\dots,e_n}$, all of the covectors $\lrc{ e_i, \cdot }$, $i\in I_{\mathrm{uf}}$, lie in a strictly convex cone.}}

\begin{proof}
This is implied by the existence of a unimodular $p^*$ map.  In Section~\ref{subsec:WXip}, I construct a particular $p^*$ map and prove its unimodularity in Proposition~\ref{prop:unimodular}.
\end{proof}

\subsection{\label{subsec:MGS}Maximal green sequence}
Let's first review what maximal green sequences are.
Recall that the $\cAp$ construction involves a ``doubled'' quiver,
where a new frozen vertex $w_i$ is introduced for each vertex $v_i$ of the original quiver for $\cA$,
and for each unfrozen vertex $v_i$ we introduce an arrow $v_i \to w_i$.
See \cite[Construction~2.11]{GHK_birational} for a more complete discussion.\footnote{The description I am giving here is for the skew-symmetric case, so at first glance it could look different from the more general construction in \cite{GHK_birational}.}
This quiver is called the {\it{framed quiver}} $\widehat{Q}$ associated to $Q$ in \cite{BDP}.
As an example, if we take principal coefficients at $s_0$ for $\ctag$,
we would replace the quiver $Q_{s_0}$ of Figure~\ref{fig:InitialQuiver} with 

\input{PrinQuiv.tex}

Now let $Q'$ be an arbitrary quiver mutation-equivalent to $\widehat{Q}$. 
An unfrozen vertex $v_i$ of $Q'$ is said to be {\it{green}} if all arrows between $v_i$ and any $w_j$ are outgoing: $v_i \rightarrow w_j$.\footnote{%
Note that by construction all unfrozen vertices of $\widehat{Q}$ are green.}
On the other hand, $v_i$ is {\it{red}} if all arrows between $v_i$ and any $w_j$ are incoming: $v_i \leftarrow w_j$.
A major result in cluster theory-- known as {\it{sign coherence of $\vb{c}$-vectors}} \cite[Corollary~5.5]{GHKK}-- implies that all unfrozen vertices are either green or red.\footnote{In the skew-symmetric case, the $\cv$-vector $\cv_{i;s}$ at the seed $s$ associated to $Q'$ is $ \cv_{i;s}:= \displaystyle{\sum_{v_i \rightarrow w_j} e_j -  \sum_{v_i \leftarrow w_j} e_j} $.  The statement that $\cv_{i;s}$ is sign coherent means either the left or right sum is empty.}
A sequence of mutations is called a {\it{green sequence}} if each mutation in the sequence is mutation at a green vertex.
It is a {\it{maximal green sequence}} if every unfrozen vertex in the resulting quiver is red.

Let $\triangle_r$ be the top $r$ rows of unfrozen vertices in $Q_{s_0}$,
and let $\vb{i}_{\triangle_r}$ be mutation at each of these vertices in order-- left to right, top to bottom.
For example, for $n=6$, $\vb{i}_{\triangle_3}$ is the following sequence: 

\input{TriangleMutation.tex}

Note that $r$ can be at most $n-2$ (the number of rows of unfrozen vertices).

{\prop{\label{prop:MGS}The sequence $\vb{i}_{\triangle_{n-2}}$, followed by $\vb{i}_{\triangle_{n-3}}, \vb{i}_{\triangle_{n-4}}, \dots, \vb{i}_{\triangle_1}$ is a maximal green sequence.}} 

{\remark{This maximal green sequence induces a simple involution on $\ssO\lrp{\cta}$ that I think is worth mentioning.  It sends $\lrp{\Vabc}^G$ to 
$\lrp{V_\alpha^* \otimes V_\gamma^* \otimes V_\beta^*}^G$.
See Corollary~\ref{cor:MGSDualizes} for details.}}

{\remark{This sequence is one of the pieces I have lifted from the earlier paper \cite{Magee}.  It has been brought to my attention that the sequence has since appeared in \cite[Section~3.3]{Le} and \cite[Proposition~9.8]{GShen_DT}.  Le discusses the sequence in the context of a {\it{cactus transformation}}, making observations similar to Corollary~\ref{cor:MGSDualizes}. Goncharov and Shen discuss the sequence in the context of Donaldson-Thomas transformations, and use it to show that an involution on $\Conf_3\lrp{\SL_n/U}$ is a cluster transformation.}}

Let's start by looking at $\itri{r}$.
Define $\trir{r}':= \lrc{w_i \left| v_i \in \trir{r}\right.}$ and 
let $F$ be the frozen vertices of $Q_{s_0}$. 
We'll split up the effects of $\itri{r}$ into three parts--
\begin{enumerate}
	\item how it affects the full subquiver with vertices $\triangle_{n-2}$,
	\item how it affects the collection of arrows between $F$ and $\triangle_{n-2}$, and
	\item how it affects the collection of arrows between $\trir{n-2}'$ and $\triangle_{n-2}$.
\end{enumerate}
Note that we {\emph{can}} split up the analysis this way.  Since we never mutate at frozen vertices, arrows between the vertices of $\trir{n-2}$ are unaffected by the presence of the frozen vertices.
There is never a composition with the center vertex frozen.
Additionally, since we never introduce arrows between frozen vertices, we could in principle treat {\emph{each}} frozen vertex separately if we wanted to.

{\lemma{\label{lem:itritri}The mutation sequence $\itri{r}$ sends the subquiver 

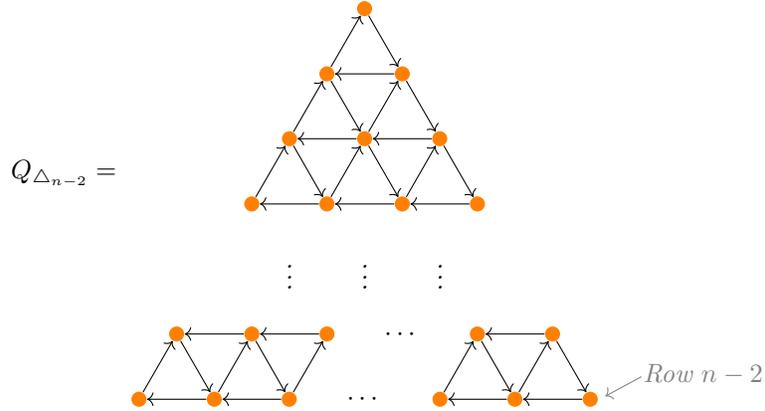
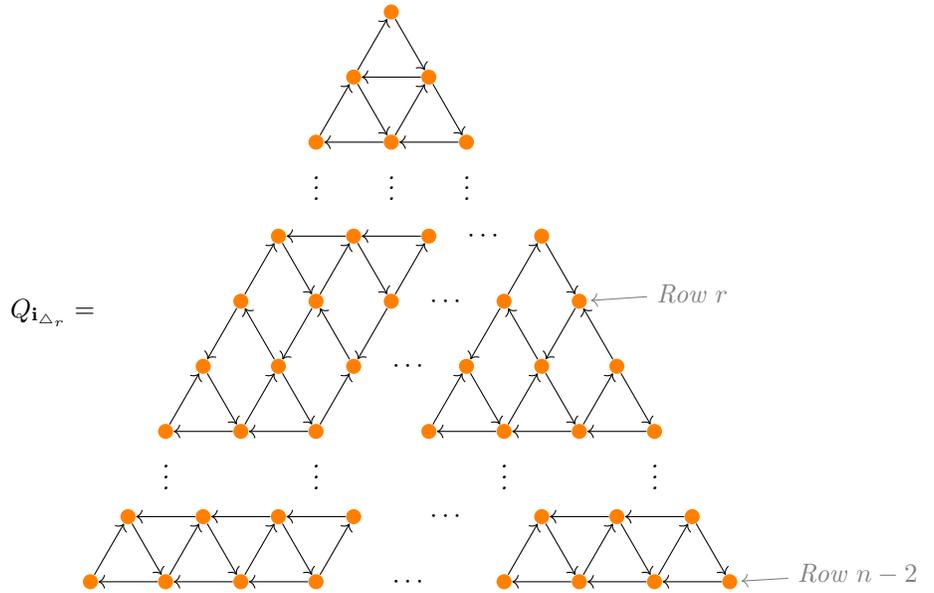
\input{Quf.tex}

to 

\input{Qufr.tex}
So $Q_{\trir{r-1}}$-- the full subquiver with vertex set $\trir{r-1}$-- remains unchanged, 
$Q_{\trir{r}}$ only has its bottom horizontal arrows deleted,
$Q_{\trir{r+1}}$ additionally has its bottom horizontal arrows deleted and its bottom diagonal arrows reversed,
and this accounts for all changes to $Q_{\trir{n-2}}$.
}}

\begin{proof}
It is immediate that the claim holds for $r=1$-- there is only one mutation to perform.
Suppose it holds for all $q < r$.
Then the quiver after performing $\itri{r-1}$ is $Q_{\itri{r-1}}$.
All that remains is mutation through row $r$.
We start with the leftmost vertex:

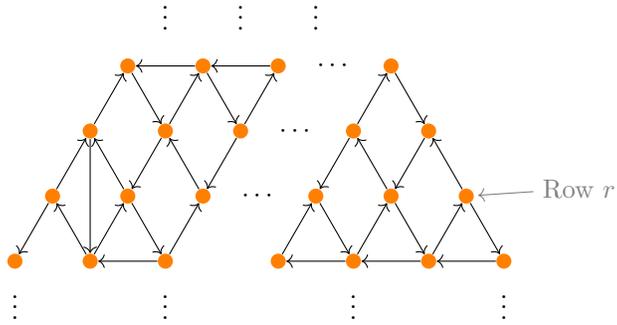
\input{Qufr1.tex}

This is followed by

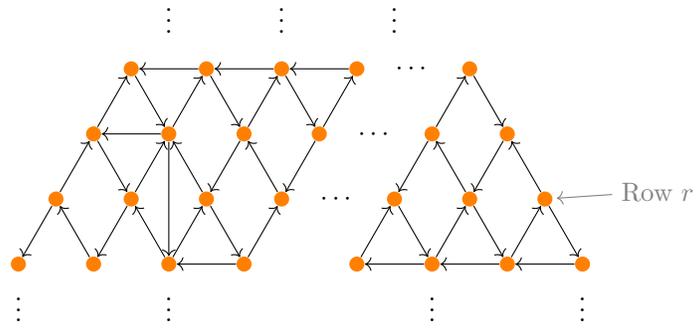
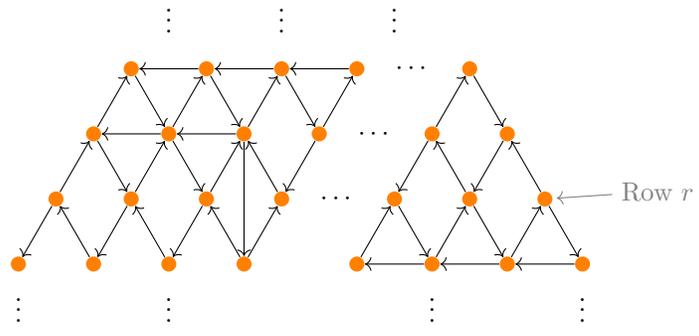
\input{Qufr2.tex}
\input{Qufr3.tex}

and so forth.
Mutation at the $k^{\text{th}}$ vertex $v_k$ of row $r$, $1 < k < r$, sends

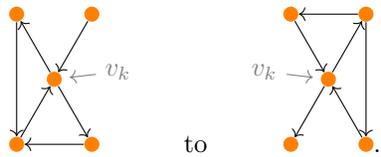
\input{Qkr.tex}
\quad to \quad
\input{Qkrmut.tex}

After mutation at $v_{r-1}$, we have the quiver
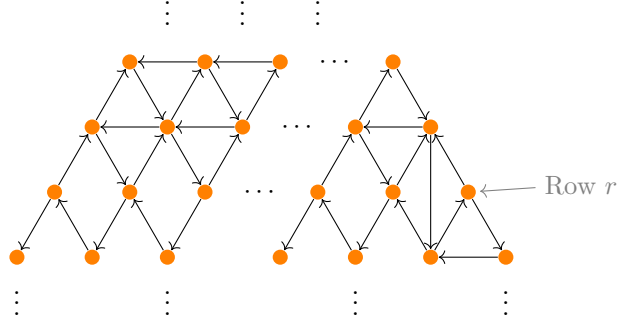
\input{Qufrr1.tex}

Finally, mutation at $v_r$ yields $Q_{\itri{r}}$.
\end{proof}

Now let's move on to how $\itri{r}$ affects the collection of arrows between $F$ and $\trir{n-2}$.
This isn't really necessary to prove Proposition~\ref{prop:MGS}-- these arrows aren't considered when determining if an unfrozen vertex is red or green--
but it's worth knowing in any case, and it will provide a nice sanity check later. See Remark~\ref{rem:reverse}.

{\lemma{\label{lem:itrifrozen}For each $Q_{s'}$ mutation equivalent to $Q_{s_0}$, let $A_{s'}$ be the subquiver having all vertices of $Q_{s'}$ but only those arrows for which either the head or tail is in $F$.
Then $A_{\itri{r}\lrp{s_0}} \subset Q_{\itri{r}\lrp{s_0}}$ can be constructed from $A_{s_0}$ as follows:
\begin{enumerate}
	\item Rearrange frozen vertices, keeping arrows fixed to their original positions. (Vertices are being relabeled.) Send 
		$v_{\lrp{n-1,1,0}}$ to the $v_{\lrp{n-r-1,0,r+1}}$ position, 
		$v_{\lrp{n-1,0,1}}$ to the $v_{\lrp{n-r-1,r+1,0}}$ position,
		$v_{\lrp{a,b,0}}$ to the $v_{\lrp{a+1,b-1,0}}$ position for $1<b<r+1$, and
		$v_{\lrp{a,0,c}}$ to the $v_{\lrp{a+1,0,c-1}}$ position for $1<c<r+1$.
	\item  Reverse arrows involving the vertices now in the 
		$v_{\lrp{n-r-1,r+1,0}}$ and $v_{\lrp{n-r-1,0,r+1}}$ positions. 
\end{enumerate}	
}}

I'll illustrate the claim with an example before proving it.
For $n=6$, $A_{s_0}$ is

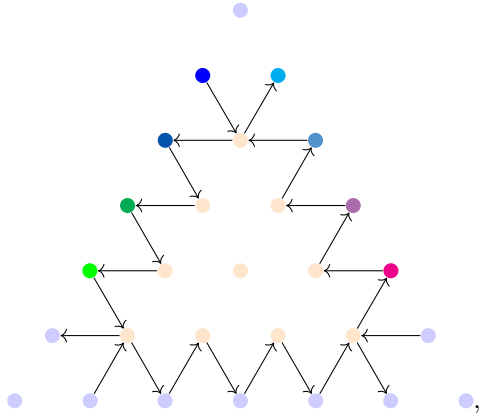
\input{As0.tex}

and $A_{\itri{3}\lrp{s_0}}$ is 

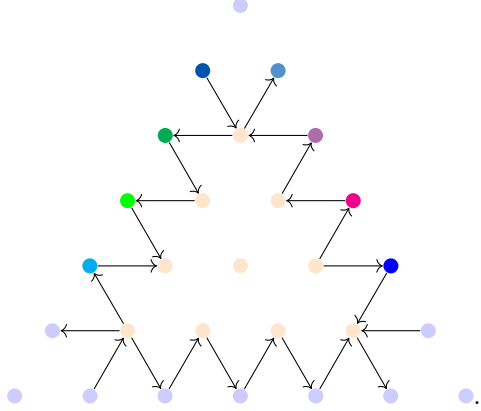
\input{Aitri.tex}

\begin{proof}
Mutating $Q_{s_0}$ at the top unfrozen vertex $v_{\lrp{n-2,1,1}}$ produces

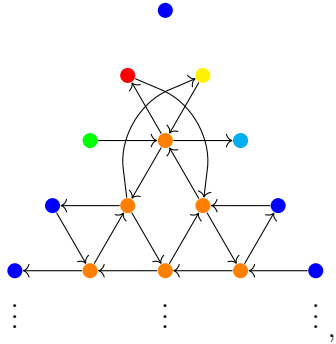
\input{Qs1.tex}
which we can rearrange as
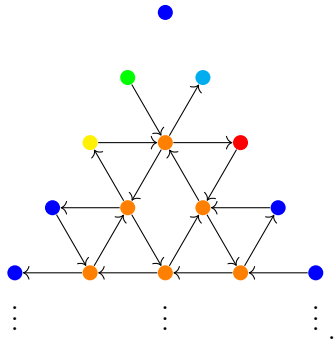
\input{Qs1rearranged.tex}
So the claim holds for $A_{\itri{1}\lrp{s_0}}$.
Suppose it holds for all $q<r$.
Then, using Lemma~\ref{lem:itritri}, 
$Q_{\itri{r-1}\lrp{s_0}}$ is
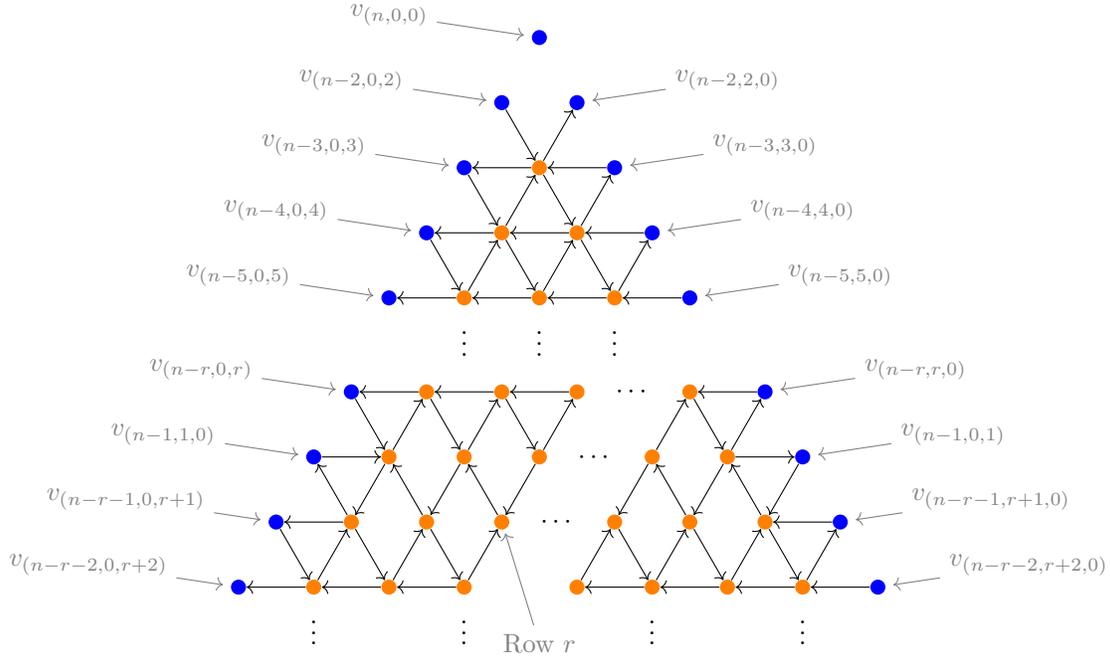
\input{Qitrir1fr.tex}
It remains to mutate through row $r$.
Mutating at the leftmost unfrozen vertex of row $r$ gives
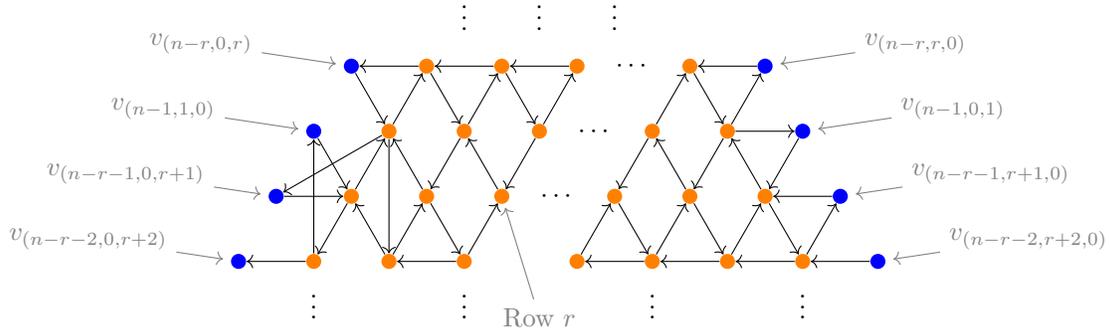
\input{Qitrir1fr1.tex}
which we rearrange as 
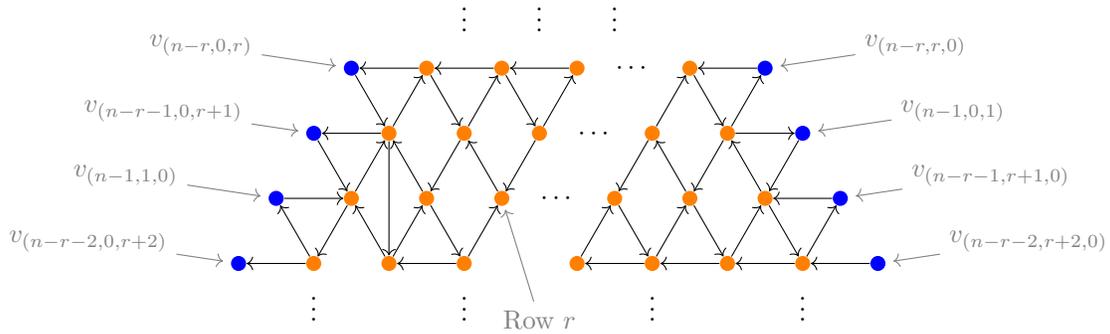
\input{Qitrir1fr1re.tex}
The next vertex of mutation sharing an arrow with some frozen vertex is the rightmost vertex of row $r$-- the final vertex in our sequence.
The penultimate quiver in the sequence is 
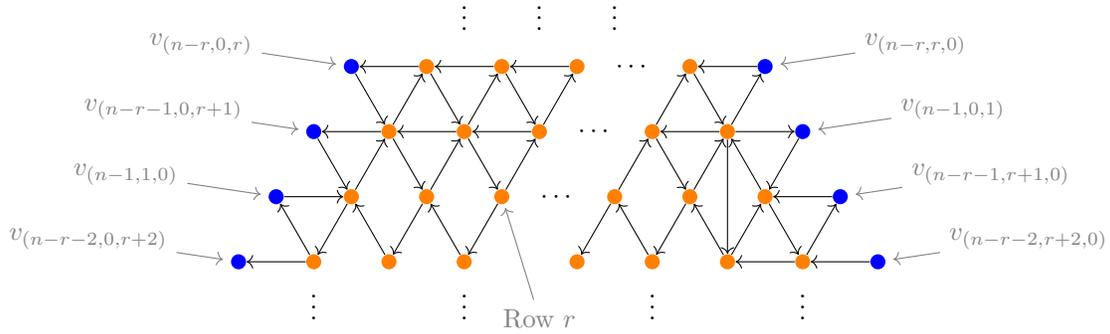
\input{Qitrir1frr1.tex}
The final mutation gives
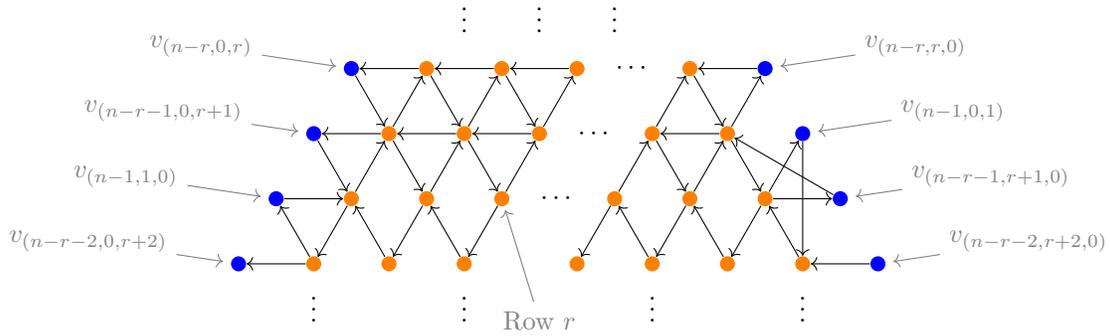
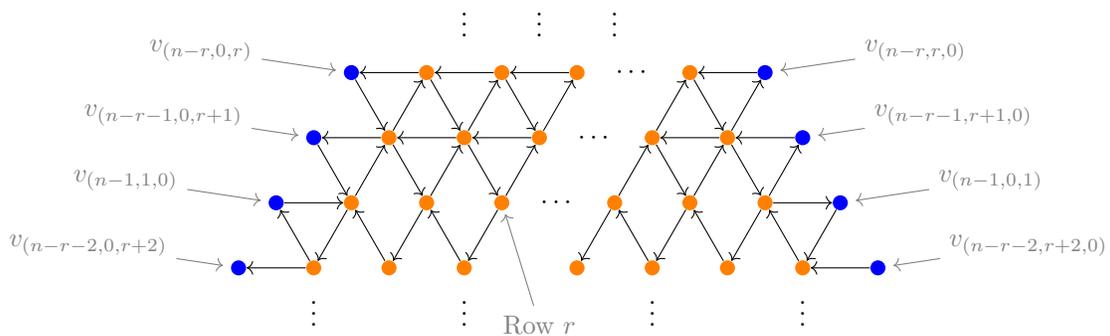
\input{Qitrir1frr.tex}
which rearranges to
\input{Qitrir1frrre.tex}
completing the proof.
\end{proof}

Now onto the arrows between $\trir{n-2}'$ and $\trir{n-2}$.

{\lemma{\label{lem:itriframed}For each $Q_{s'}$ mutation equivalent to $\widehat{Q_{s_0}}$, let $R_{s'}$ be the subquiver with vertex set $\trir{n-2} \bigcup \trir{n-2}'$ but only those arrows for which either the head or tail is in $\trir{n-2}'$.
Then $R_{\itri{r}\lrp{s_0}} $ can be constructed from $R_{s_0}$ as follows:
\begin{enumerate}
	\item Rearrange $\trir{n-2}$, keeping arrows fixed to their original positions. (Vertices are being relabeled.) Send 
		$w_{\lrp{n-b-1,b,1}}$ to the $w_{\lrp{n-r-1,b,r-b+1}}$ position for $b \leq r$ 
		and
		$w_{\lrp{a,b,n-a-b}}$ to the $w_{\lrp{a+1,b,n-a-b-1}}$ position for $b \leq r$. 
	\item  Reverse the arrow between $v_{\lrp{n-r-1,b,r-b+1}}$ and the vertex now in the $w_{\lrp{n-r-1,b,r-b+1}}$ position for $b \leq r$.
	\item  Introduce a new arrow from the vertex now in the $w_{\lrp{a+1,b,n-a-b-1}}$ position to $v_{\lrp{n-r-1,b,r-b+1}}$ for $b \leq r$.
	\item  If $r<n-2$, introduce a new arrow from $v_{\lrp{n-r-2,b,r-b+2}}$ to the vertex now in the $w_{\lrp{a,b,n-a-b}}$ position for $b \leq r$.
\end{enumerate}	
}}
To illustrate the claim, if we take $n=7$, then $R_{s_0}$ is

\input{Qframed.tex}

and $R_{\itri{3}\lrp{s_0}} $ is

\input{Qframedprime.tex}

\begin{proof}
The first mutation gives

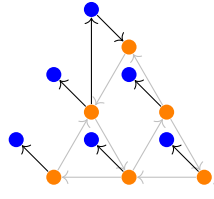
\input{Qframut1.tex}

which agrees with the statement for $r=1$.
Assume it holds for all $q<r$.
Then after mutating through $\itri{r-1}$ and rearranging the vertices as described, we have

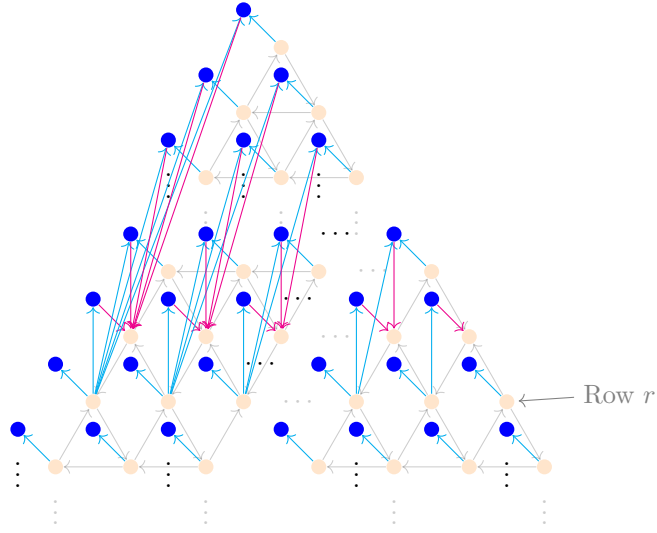
\input{Qframutr1.tex}

and we just have to mutate through row $r$.
The unfrozen portion of the quiver for each of the remaining mutations is given in the proof of Lemma~\ref{lem:itritri}.
Note that there is a {\textcolor{cyan}{cyan}} arrow emanating from $v_{\lrp{n-r-1,b,r-b+1}}$ corresponding to each 
{\textcolor{magenta}{magenta}} arrow terminating at $v_{\lrp{n-r-2,b,r-b+2}}$,
and there is one additional \tc{cyan}{cyan} arrow $v_{\lrp{n-r-1,b,r-b+1}} \tc{cyan}{\rightarrow} w_{\lrp{n-r-1,b,r-b+1}} $.
Now, $v_{\lrp{n-r-1,b,r-b+1}}$ is the $b^{\text{th}}$ vertex of mutation in this row, and each of the 
{\textcolor{magenta}{magenta}} arrows are killed by a composition 
$$v_{\lrp{n-r-2,b,r-b+2}} \rightarrow v_{\lrp{n-r-1,b,r-b+1}} \tc{cyan}{\rightarrow} \tc{blue}{\bullet}$$
while a \tc{cyan}{cyan} arrow $v_{\lrp{n-r-2,b,r-b+2}} \tc{cyan}{\rightarrow} w_{\lrp{n-r-1,b,r-b+1}} $ is created.
Meanwhile, if $r<n-2$, a new \tc{cyan}{cyan} arrow is created by the compositions
\eqn{
v_{\lrp{n-r,b,r-b}} \rightarrow v_{\lrp{n-r-1,b,r-b+1}} \tc{cyan}{\rightarrow} \tc{blue}{\bullet} \qquad \RA \qquad v_{\lrp{n-r,b,r-b}} \tc{cyan}{\rightarrow} \tc{blue}{\bullet},
}
and each of the \tc{cyan}{cyan} arrows
\eqn{
v_{\lrp{n-r-1,b,r-b+1}} \tc{cyan}{\rightarrow} \tc{blue}{\bullet}
}
is reversed, becoming the \tc{magenta}{magenta} arrow
\eqn{
v_{\lrp{n-r-1,b,r-b+1}} \tc{magenta}{\leftarrow} \tc{blue}{\bullet}.
}
So after performing $\itri{r}$ we obtain the quiver

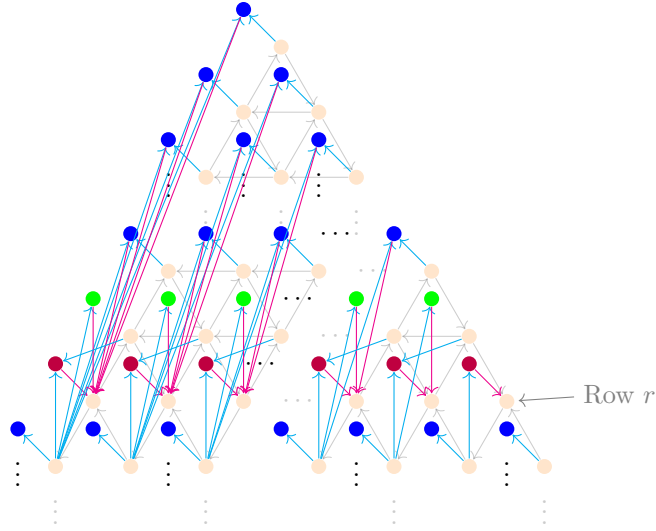
\input{Qframutr.tex}

which we rearrange to 

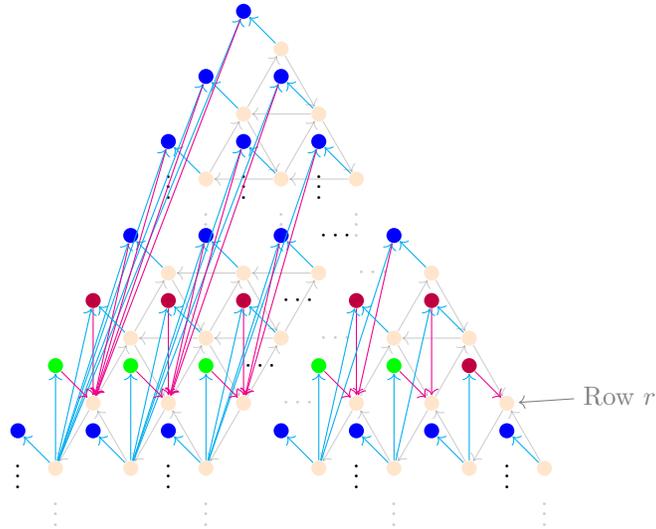
\input{Qframutrre.tex}

finishing the proof.
\end{proof}

We now have all of the ingredients we need to tackle Proposition~\ref{prop:MGS}.
\begin{proof}
We start with $\itri{n-2}$.
From the proof of Lemma~\ref{lem:itriframed},
we see that each time we mutate at a vertex $v_k$ in $\itri{n-2}$, all arrows between $v_k$ and $\trir{n-2}'$ emanate from $v_k$-- so $v_k$ is green.
Then $\itri{n-2}$ is a green sequence.
Using Lemmas~\ref{lem:itritri}, \ref{lem:itrifrozen}, and \ref{lem:itriframed},
performing $\itri{n-2}$ on $\widehat{Q_{s_0}}$ and rearranging vertices as indicated in the lemmas results in the quiver

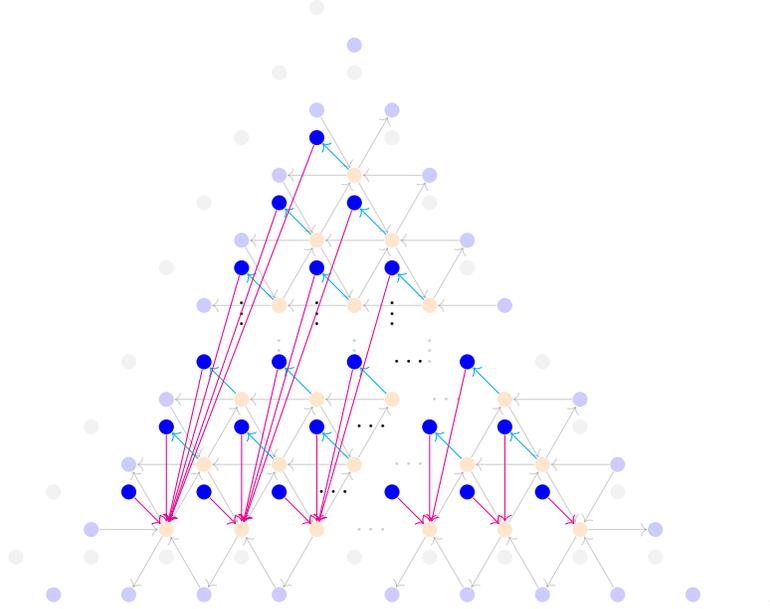
\input{Qhatn2.tex}

The vertices of the bottom unfrozen row are now red, while the remaining unfrozen vertices are all green.
For consistency with the $\trir{r}$ notation, let's only consider unfrozen vertices when indexing the rows.
So the bottom unfrozen row we'll call row $n-2$, the one above it row $n-3$, and so forth.
This quiver is very similar to the one we started with.
Above row $n-2$ the only relevant difference is the introduction of the \tc{magenta}{magenta} arrows from $\trir{n-3}'$ to row $n-2$.
Referring to the proof of Lemma~\ref{lem:itritri}, we note that no vertex of mutation in the sequence $\itri{n-4}$ shares an arrow with row $n-2$.
As a result, no composition affecting these arrows can occur until we mutate at row $n-3$.
That is, the subsequence $\itri{n-4}$ of $\itri{n-3}$ proceeds exactly as before, with these \tc{magenta}{magenta} arrows tagging along for the ride.
Then prior to mutation through row $n-3$, there is a \tc{cyan}{cyan} arrow emanating from $v_{\lrp{2,b,n-b-2}}$ for all but one of the \tc{magenta}{magenta} arrows
terminating at $v_{\lrp{1,b,n-b-1}}$.
These paired \tc{magenta}{magenta} arrows are canceled upon mutation at  $v_{\lrp{2,b,n-b-2}}$.
So after performing $\itri{n-3}$ and rearranging frozen vertices, we have

\input{Qhatn3.tex}

Now the unfrozen vertices of rows $n-2$ and $n-3$ are red and the remaining unfrozen vertices are green.
We can employ the reasoning just used for $\itri{n-3}$ to the remaining subsequences $\itri{n-4}, \itri{n-5}, \dots, \itri{1}$. 
The resulting quiver is

\input{Qhatfinal.tex}

and the sequence $\itri{n-2}$, followed by $\itri{n-3},\itri{n-4},\dots,\itri{1}$ is a maximal green sequence.
\end{proof}

{\remark{\label{rem:reverse}With the indicated rearranging of frozen vertices, the final quiver we obtained is the same as the original framed quiver with every arrow reversed.
Imagine each $w_{\lrp{a,b,c}}$ as lying above $v_{\lrp{a,b,c}}$.
Now ignore temporarily all vertices that aren't attached to any arrows, 
and reflect the rest of the final quiver over the plane given by $a=c$.
The quiver itself is obviously the same.
We've just changed its embedding into $\R^3$
and returned each $w_{\lrp{a,b,c}}$ to its original position.
Note that this also gives an isomorphism of the full subquiver whose vertex set is all of the $v_{\lrp{a,b,c}}$'s with the quiver $Q_{s_0}$.
So there is an isomorphism of the final quiver with the {\it{coframed quiver}}\footnote{This differs from the framed quiver in that arrows $w_i \to v_i$ are introduced rather than $v_i\to w_i$.} $\widecheck{Q_{s_0}}$ fixing the $w_{\lrp{a,b,c}}$'s.
This is what we expect by \cite[Proposition~2.10]{BDP}, and it is a sanity check for the work in this section. 
}}
%

\section{\label{sec:compactify}From \texorpdfstring{$\ctag$}{generic configurations} to \texorpdfstring{$\cta$}{configurations}}
\subsection{\label{subsec:optimized}Existence of optimized seeds}
The main result of this subsection is that every frozen index for $\ctag$ has an {\it{optimized seed}}.
In the skew-symmetric case, a seed $s$ is {\it{optimized}} for the frozen index $i$ if the vertex $v_i$ is a sink of the associated quiver $Q_s$.
The frozen index $i$ corresponds to a component $D_i$ of the divisor we've added to $\ctag$,
and the existence of a seed optimized for $i$ ensures that the $\vartheta$-function on $\ctag^\vee$ associated to $D_i$ is a global monomial.
See \cite[Section~9.1]{GHKK} for a more complete discussion of optimized seeds.

So far we have a basis $\vb{B}^\times$ for $\ssO\lrp{\ctag}$.
What we really want is a basis $\vb{B}$ for $\ssO\lrp{\cta}$.
A natural candidate for $\vb{B}$ is the subset of $\vb{B}^\times$ that extends to the divisors we've added, {\it{i.e.}} $\vb{B}^\times \bigcap \ssO\lrp{\cta}$.
But this candidate isn't automatically a basis for $\ssO\lrp{\cta}$.
Maybe $\vartheta_p, \vartheta_q \in \vb{B}^\times$ both have a pole along some component $D_i$, but these poles cancel in their sum $\vartheta_p+\vartheta_q$.
Then we would have $\vartheta_p+\vartheta_q \in \ssO\lrp{\cta}$, but $\vartheta_p,\vartheta_q \notin \ssO\lrp{\cta} \subset \ssO\lrp{\ctag}$.
The existence of an optimized seed for each frozen index ensures that this does not happen-- if a linear combination of $\vartheta$-functions extends to $D_i$, 
then each $\vartheta$-function in the sum extends as well.\cite[Proposition~9.7]{GHKK}
This condition is needed to utilize \cite[Theorem~0.19]{GHKK}, which will be used in the coming subsection on the potential $W$ and cone $\Xi$ for $\cta$.

{\prop{\label{prop:op}Every frozen index for $\ctag$ has an optimized seed.}}

\begin{proof}
For cluster varieties with skew-symmetric exchange matrix, a seed $s$ is optimized for the frozen index $f$ if and only if the vertex $v_f$ is a sink in the quiver $Q_s$.\cite[Lemma~9.2]{GHKK}
Consider a quiver $Q_L$ of the form
\input{QLine.tex}
The sequence of mutations $v_1,v_2, \dots, v_r$ yields the quiver
\input{QLineMut.tex}
making $v_f$ a sink.
The initial seed quiver for $\ctag$ is shown in Figure~\ref{fig:InitialQuiver}.
Call it $Q_{s_0}$.
Since there are no arrows to or from the corner vertices $v_{\lrp{n,0,0}}$, $v_{\lrp{0,n,0}}$, and $v_{\lrp{0,0,n}}$,
every quiver mutation equivalent to $Q_{s_0}$ will trivially be optimized for these three vertices.
Beyond that, $Q_{s_0}$ is optimized for $v_{\lrp{n-1,1,0}}$, $v_{\lrp{0,n-1,1}}$, and $v_{\lrp{1,0,n-1}}$.
For the remaining frozen vertices $v_f$, there is a subquiver of $Q_{s_0}$ isomorphic to $Q_L$.
Performing these mutations on $Q_{s_0}$ only affects the subquiver whose vertices are either in $Q_L$ or connected to $Q_L$ by an arrow.
As arrows between frozen vertices are deleted,
any frozen vertices besides $v_f$ can be ignored when determining if $v_f$ becomes a sink.
Then the relevant subquiver of $Q_{s_0}$ has the form 

\input{QSub.tex}

possibly with the top or bottom row deleted and with the middle row being the subquiver $Q_L$.
(Of course, depending upon the position of $v_f$, it may be necessary to rotate Figure~\ref{fig:SubQuiv}.)
The key observation is that, for every quiver in the sequence, each vertex connected to $v_f$ by an arrow is a vertex of the subquiver $Q_L$.
The cycles prevent any new arrows involving $v_{f}$ from developing via some composition with an arrow not in $Q_L$.
The explicit mutations, ending with $v_f$ as a sink, are shown below.
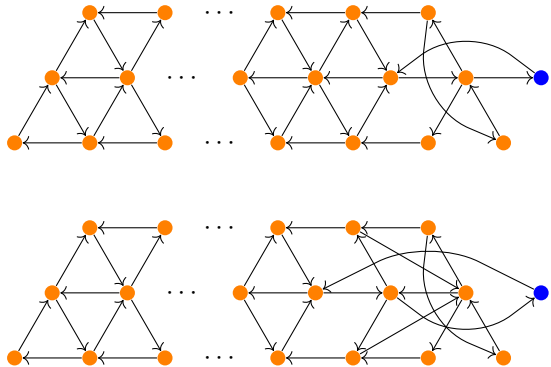
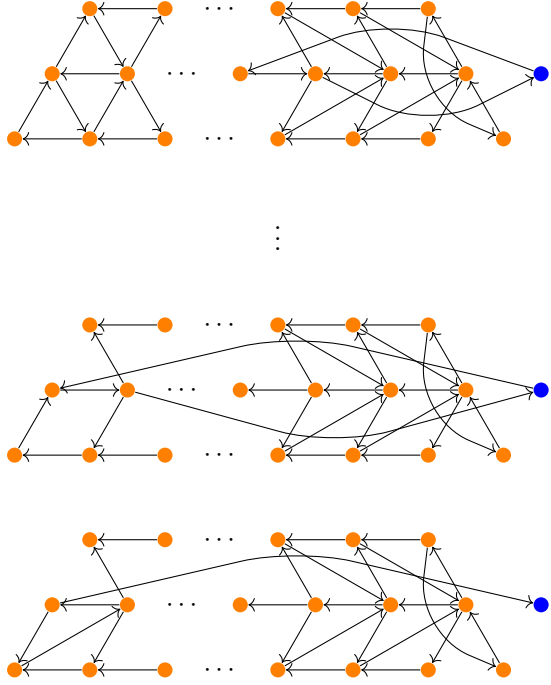
\input{QSubmut1.tex}
\input{QSubmut2.tex}
\input{QSubmut3.tex}
\eqn{\vdots}
\input{QSubmutn-1.tex}
\input{QSubmutn.tex}
\end{proof}

\subsection{The potential \texorpdfstring{$W$}{W} and cone \texorpdfstring{$\Xi$}{Xi} for \texorpdfstring{$\cta$}{configurations of flags}}
In this subsection I compute the Landau-Ginzburg potential $W$ and corresponding cone ${\Xi := \lrc{W^T \geq 0} \subset \ctag^\vee\lrp{\R^T}}$.
By \cite[Theorem~0.19]{GHKK}, the analogous cone for $\overline{\cA}_\prin$ gives a canonical basis for the finitely generated algebra $\cmid\lrp{\overline{\cA}_\prin}=\up\lrp{\overline{\cA}_\prin}$.
However, the exchange matrix for $\ctag$ is full rank over $\Z$.  This is immediate from the stronger result Proposition~\ref{prop:unimodular}.
Then, as explained in \cite[Proofs of Corollaries~0.20 and 0.21, page~581]{GHKK}, the desired results for $\cta$ are implied by the results for $\overline{\cA}_\prin$.
I'll say a few more words about this in Subsection~\ref{subsec:Tk}.

I give an explicit description of $W$ and $\Xi$ in the initial seed.
In this seed, the inequalities defining $\Xi$ are precisely the tail positivity conditions of \cite{Zel_Tail}.
I also exhibit a map $p^*:N \to M$ that identifies $W$ with 
the representation theoretically defined potential $\mathcal{W}_{GS}$ of \cite{GShen} on $\ctag$ and identifies $\Xi$ with the Knutson-Tao hive cone.\footnote{Here $N$ is the cocharacter lattice for a torus in $\ctag$ and $M$ is its dual-- the cocharacter lattice for a torus in $\ctag^\vee$.  The map $p^*$ commutes with mutation and is closely related to the exchange matrix.  See \cite[Section~2]{GHK_birational} for a general discussion of $p^*$ maps, and Subsection~\ref{subsec:WXip} below for a brief review.}\cite{KTv1}
Before doing this, let's recall what $\mathcal{W}_{GS}$ and the Knutson-Tao hive cone are.
\subsubsection{Knutson-Tao hive cone}

Consider a triangular array of vertices indexed by triples $\lrp{a,b,c} \in \lrp{\Z_{\geq 0}}^3$ with $a + b + c =n$, for some fixed $n$, just like in Figure~\ref{fig:confIndexed}.
Let $\mathcal{H}$ be the set of these vertices.
$\R^{\mathcal{H}}$ is the possible labelings of these vertices by real numbers.
Now take any pair of neighboring triangles, together forming a rhombus.
This rhombus defines a linear inequality in $\R^{\mathcal{H}}$ by requiring
the sum of the labels on the obtuse vertices to be greater than or equal to the sum of the labels on the acute vertices.

\input{RhombusIneq.tex}

Now denote by ${\bf{1}}_{\mathcal{H}}$ the labeling where each entry is 1.
The ``Knutson-Tao hive cone'' generally refers to one of the following three cones:
\begin{enumerate}
        \item \label{KTbig} the polyhedral cone in $\R^{\mathcal{H}}$ satisfying all rhombus inequalities
        \item \label{KTslice} the slice of (\ref{KTbig}) having top entry 0
        \item \label{KTquot} the quotient of (\ref{KTbig}) by $\R \cdot {\bf{1}}_{\mathcal{H}}$.\footnote{Note that ${\bf{1}}_{\mathcal{H}}$ spans a linear subspace of (\ref{KTbig}).}
\end{enumerate}
(\ref{KTslice}) and (\ref{KTquot})
clearly have completely equivalent combinatorics, with each point in
(\ref{KTslice})
giving a representative of one of the equivalence classes in (\ref{KTquot}).
The points in the Knutson-Tao hive cone are called {\it{hives}}.

The Knutson-Tao hive cone encodes the Littlewood-Richardson coefficients in a really beautiful way.
Suppose we want to know
$\dim{\lrp{V_\alpha \otimes V_\beta \otimes V_\gamma}^G}$.
The choice of weights $\lrp{\alpha, \beta, \gamma }$ determines the border of a hive, which I'll illustrate in terms of (\ref{KTquot}).
If $\lambda = \lrp{\lambda_1, \dots, \lambda_n}$, define $\lrm{\lambda}= \lambda_1 + \cdots + \lambda_n$.
Now take $\lrm{\alpha} + \lrm{\beta} + \lrm{\gamma} = 0$-- otherwise
$\dim{\lrp{V_\alpha \otimes V_\beta \otimes V_\gamma}^G}=0$.
Then we label the border of the hive as follows:

\input{HiveBorder.tex}

Note that the condition $\lrm{\alpha} + \lrm{\beta} + \lrm{\gamma} = 0$ is exactly what we need to be able to fill in the border this way.
Also note that, since we are working only up to translations in the ${\bf{1}}_{\mathcal{H}}$ direction, the border is completely determined by the choice of $\lrp{\alpha, \beta, \gamma}$.
Furthermore, this picture is manifestly symmetric under cyclically permuting $\alpha$, $\beta$, and $\gamma$.
Knutson and Tao showed that the number of integral hives with this border is precisely
$\dim{\lrp{V_\alpha \otimes V_\beta \otimes V_\gamma}^G}$.\cite{KTv1,Buch}

\subsubsection{Goncharov-Shen potential \texorpdfstring{$\mathcal{W}_{GS}$}{WGS}}

In \cite{GShen}, Goncharov and Shen gave a new construction of the Knutson-Tao hive cone, which I describe briefly here.
Goncharov and Shen describe points in $\Fld$ as pairs $\lrp{U,\chi}$, where $U$ is a
maximal unipotent subgroup in $G$ and $\chi$ is
a non-degenerate additive character on $U$-- meaning a group homomorphism ${U \to \C_a}$
such that the stabilizer of $\lrp{U, \chi}$ under the conjugation action of $G$ is precisely $U$.
For each triple $\lrp{\lrp{U_1 , \chi_1 }, \lrp{U_2 , \chi_2 }, \lrp{U_3 , \chi_3 }} \in \ctag$,
there is a unique element $u_{jk} \in U_i$ conjugating $U_j$ to $U_k$.
This gives a natural function on $\ctag$:
\eqn{
\mathcal{W}_{GS} \lrp{\lrp{U_1 , \chi_1 }, \lrp{U_2 , \chi_2 }, \lrp{U_3 , \chi_3 }} := \chi_1 \lrp{u_{23}} +  \chi_2 \lrp{u_{31}} +  \chi_3 \lrp{u_{12}}.
}
They then show that in the initial seed of the cluster variety,
$\mathcal{W}_{GS}^T \geq 0$ gives exactly the rhombus inequalities cutting out the Knutson-Tao hive cone.

\subsubsection{\label{subsec:WXip}\texorpdfstring{$W$}{W}, \texorpdfstring{$\Xi$}{Xi}, and \texorpdfstring{$p^*$}{p}}

The Landau-Ginzburg potential $W$ is the sum of $\vartheta$-functions associated to the irreducible components of $D:= \cta \setminus \ctag$.
$D$ is given by
\eqn{
\sum_{i=1}^{n-1} \lrp{ D_{\lrp{i,n-i,0}} + D_{\lrp{0,i,n-i}} + D_{\lrp{n-i,0,i}} }, 
}
where, for example, 
\eqn{
D_{\lrp{i,n-i,0}} :=& \lrc{\lrp{X,Y,Z} \in \cta \left| X_i \not\pitchfork Y_{n-i} \right.} \footnotemark\\
		  =&\lrc{ A_{\lrp{i,n-i,0}} = 0 }.
}
\footnotetext{The symbol $\not\pitchfork$ denotes a non-transverse intersection.}

Suppose the seed $s$ is optimized for the frozen index $\lrp{i,n-i,0}$.
Then on the torus $T_{M;s}$ in the atlas for $\ctag^\vee$, 
$\vartheta_{\lrp{i,n-i,0}}$ is given by $z^{-e_{\lrp{i,n-i,0}}}$.\cite[Lemma~9.3]{GHKK}\footnote{The negative sign comes from the sign change identification $i$ of 
$\ctag^\vee\lrp{\Z^T}$ and $\lrp{\ctag^\vee}^\trop\lrp{\Z}$.}
We can express $\vartheta_{\lrp{i,n-i,0}}$ on other tori in the atlas by pulling back $z^{-e_{\lrp{i,n-i,0}}}$ via the birational gluing maps.
The formula for mutation at $v_k$ is
\eqn{
\mu_k^*\lrp{ z^{\vb{n}}} = z^{\vb{n}} \lrp{ 1 + z^{e_k}}^{-\lrc{\vb{n},e_k}},
}
where $\vb{n} \in N$-- the lattice of the fixed data used to define the cluster structure.
If $s= \lrp{e_1, \dots, e_n}$, then $\mu_k \lrp{s} = \lrp{e_1',\dots, e_n'}$ where
\eqn{
e_i' = \begin{cases}
e_i + \lrb{\epsilon_{ik}}_+ e_k & {\mbox{if }} i\neq k\\
-e_k & {\mbox{if }} i=k
\end{cases}.
}
Here $\epsilon$ is the exchange matrix for the seed $s$, given by $\epsilon_{ij}= \lrc{e_i,e_j}$.\footnote{Generally there would be multipliers $d_j$ modifying this expression, but in the case under consideration all multipliers $d_j$ are $1$.  This is known as the skew-symmetric case.  Additionally, the reader who wishes to compare results here with work following Fomin and Zelevinsky's conventions should set the Fomin-Zelevisky exchange matrix $B$ equal to $\epsilon^T$.}
The notation $\lrb{a}_+$ is shorthand for $\min\lrc{0,a}$.
I'll express each of the $\vartheta$-functions, and hence $W$, in the initial seed $s_0$ using this mutation formula.

{\remark{Using this exponential notation, the $\cX$-variables of a given seed $s = \lrp{e_1, \cdots, e_n}$ are defined by $X_i := z^{e_i}$ and the $\cA$-variables by $A_i :=z^{e_i^*}$. See \cite[Section~2]{GHK_birational}, keeping in mind that we are considering the skew-symmetric case here.}}

For each frozen index $f$, we have an explicit sequence of mutations from $s_0$ to a seed $s_f$ optimized for $f$ from the proof of Proposition~\ref{prop:op}.
We want to pullback $\vartheta_f$ from $T_{M;s_f}$ to $T_{M;s_0}$, so we reverse this sequence.
All of the mutations occur at vertices of the subquiver $Q_L$,
so, by the mutation formula given above, only indices of $Q_L$ will come into play in computing the pullback of $\vartheta_f$.

{\prop{\label{prop:theta}
Recall the quivers $Q_L$ and $Q_{L_f}$ of Proposition~\ref{prop:op}.
Call the seeds associated to these quivers $s_0$ and $s$.
Then the pullback of $z^{-e_{f_{s}}}$ from $T_{M;s}$ to $T_{M;s_0}$ is 
\eqn{
z^{-e_f} +
z^{-e_f-e_1} +
z^{-e_f-e_1-e_2} +
\cdots +
z^{-e_f-e_1-e_2 -\cdots - e_r}.
}
}}

\begin{proof}
The quiver for the first mutation is
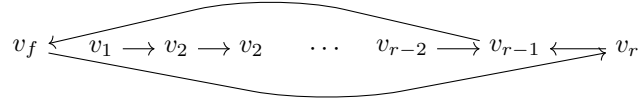
\input{QLineMut1.tex}
So
\eqn{
\mu_{r}^* \lrp{z^{-e_f'}} &= z^{-e_f'}\lrp{1+z^{e_r}}^{-\lrc{-e_f',e_r}}\\
&= z^{-e_f - \lrb{ \epsilon_{fr}}_+ e_r} \lrp{1 + z^{e_r}}^{\lrc{e_f + \lrb{\epsilon_{fr}}_+e_r,e_r}}\\
&= z^{-e_f - e_r} \lrp{1 + z^{e_r}}.
}
The next quiver is
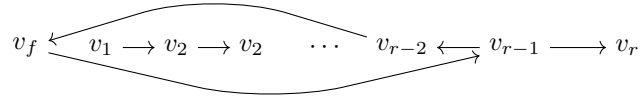
\input{QLineMut2.tex}
\eqn{
\mu_{r-1}^* \lrp{z^{-e_f' - e_r'} \lrp{1 + z^{e_r'}}} &=
z^{-e_f' - e_r'}\lrp{1+z^{e_{r-1}}}^{-\lrc{-e_f'-e_r',e_{r-1}}} \lrp{1 + z^{e_r'} \lrp{1+z^{e_{r-1}}}^{-\lrc{e_r',e_{r-1}}}}\\
&=
z^{-e_f - e_{r-1} - e_r }\lrp{1+z^{e_{r-1}}}^0 \lrp{1 + z^{e_r} \lrp{1+z^{e_{r-1}}}^{1}}\\
&=
z^{-e_f - e_{r-1} - e_r }\lrp{1 + z^{e_r} \lrp{1+z^{e_{r-1}}}}.
}
This pattern continues with the $i^{\text{th}}$ mutation yielding
\eqn{
z^{-e_f-e_{r-i+1} - e_{r-i+2} - \cdots - e_{r}} \lrp{1 + z^{e_r} \lrp{ 1 + z^{e_{r-1}} \cdots \lrp{1+z^{e_{r-i+1}}} \cdots }}.
}
The result after all $n$ mutations is
\eqn{
& z^{-e_f-e_1 - \cdots - e_r} \lrp{1 + z^{e_r} \lrp{ 1 + z^{e_{r-1}} \cdots \lrp{1+z^{e_{1}}} \cdots }}\\
=& z^{-e_f-e_1 - \cdots - e_r}
+  z^{-e_f-e_1 - \cdots - e_{r-1}}
+ \cdots + z^{-e_f},
}
as claimed.
\end{proof}

Using Proposition~\ref{prop:theta}, we can immediately express $W$ on the torus $T_{M;s_0}$ of the initial seed of $\ctag^\vee$.

{\cor{\label{cor:W}
Take $a,b,c \in \Z_{> 0}$.
The restriction of $W$ to $T_{M;s_0}$ is
\eqn{
W= \sum_{a+b = n} \vartheta_{\lrp{a,b,0}} 
+\sum_{b+c = n} \vartheta_{\lrp{0,b,c}}
+ \sum_{a+c = n} \vartheta_{\lrp{a,0,c}}, 
}
where  
\eqn{
\vartheta_{\lrp{a,b,0}}= 
\sum_{i= 0}^{n-a-1} z^{-\sum_{j=0}^i e_{\lrp{a,b-j,j}}}, 
}
\eqn{
\vartheta_{\lrp{0,b,c}}= 
\sum_{i= 0}^{n-b-1} z^{-\sum_{j=0}^i e_{\lrp{j,b,c-j}}}, 
}
and
\eqn{
\vartheta_{\lrp{a,0,c}}= 
\sum_{i= 0}^{n-c-1} z^{-\sum_{j=0}^i e_{\lrp{a-j,j,c}}}.
}
}}

Note that we now have the basis $\vb{B}$ of $\ssO\lrp{\cta}$ that we were after--
it is canonically identified with $$\Xi\lrp{\Z^T} := \lrc{W^T \geq 0} \bigcap \ctag^\vee\lrp{\Z^T}.$$
Explicitly, the $\tf$-functions $\tf_p$ is in $\vb{B}$ if and only if $p$ is in $\Xi\lrp{\Z^T}$.
This basis is canonically determined by the pair $\ctag \subset \cta$.
The subset $\ctag$ is invariant under the $H^{\times 3}$ action on $\cta$,
so $\vb{B}$ must be preserved by this action.
Since $\vb{B}$ is a discrete set and $H^{\times 3}$ acts continuously,
the only possibility is that each element of $\vb{B}$ is fixed by $H^{\times 3}$.
The elements of $\vb{B}$ are defined up to scaling, so this means that every element of $\vb{B}$ is an $H^{\times 3}$-eigenfunction.
The $H^{\times 3}$ action and the weights of basis elements under this action are discussed further in Subsection~\ref{subsec:Tk}.

At this point we'd like to see if $W$ to pulls back to the Goncharov-Shen potential $\mathcal{W}_{GS}$ on $\ctag$ for some carefully chosen $p^*$.
The guideline for writing down this map will be the representation theoretic interpretation of the cones on both sides.
For this, we'll compare version (\ref{KTquot}) of the Knutson-Tao hive cone to $\Xi$.
To have a nice representation theoretic interpretation of $\Xi$, we need to relate the {\it{$\gv$-vector}} of a $\vartheta$-function to its weight under the $H^{\times 3}$ action.

First let's review what $\gv$-vectors are and fix some notation.  
Given any cluster $\cA$-variety, 
Gross-Hacking-Keel-Kontsevich build a family of deformations of $\cA$, denoted $\cAp$, that itself has the structure of an $\cA$-variety.
I recall how to build a quiver for $\cAp$ at the beginning of Section~\ref{subsec:MGS}.
A choice of seed $s$ determines a cluster torus $T_{N;s}$ in $\cA$,  an action of $T_{N;s}$ on $\cAp$, and a canonical extension of each cluster monomial to $\cAp$ (simply by performing mutations in $\cAp$ rather than $\cA$).\footnote{I'm restricting to the skew-symmetric case here to avoid a discussion of multipliers $d_i$ and to simplify notation.}
The $T_{N;s}$ action is induced by the inclusion
\eqn{
N &\hookrightarrow N\oplus M\\
n &\mapsto \lrp{n, p^*(n)}.
}
The extended cluster monomials are all eigenfunctions of this $T_{N;s}$ action, and the weight under this action is known as the $\gv$-vector at $s$ for the cluster monomial.
The choice of seed $s$ also determines an identification of $\cX\lrp{\Z^T}$ with $M$, and this weight (an element of $M$) is how we explicitly associate a $\vartheta$-function to a point in $\cX\lrp{\Z^T}$.

If $S$ is a subset of some real tropical space $\mathcal{U}\lrp{\R^T}$, define $S\lrp{\Z^T}$ to be its $\Z^T$ points-- $S\lrp{\Z^T} := S \bigcap \mathcal{U}\lrp{\Z^T}$.
Now let $p \in \Xil \lrp{\Z^T}$ and express the $\gv$-vector of $\tf_p$ at $s_0$ as 
\eqn{
\gv_{s_0} \lrp{\vartheta_p} = \sum_{a+b+c = n} g_{\lrp{a,b,c}} e_{\lrp{a,b,c}}^*. 
}
Here $\lrc{e_{\lrp{a,b,c}}^*}_{\lrp{a,b,c} \in \mathcal{H}}$ is the dual basis to the ordered basis $s_0$ of $N$.
Now, $\gv_{s_0} \lrp{\vartheta_p}$
is the exponent of the leading term of $\vartheta_p$ expressed as a Laurent polynomial on $T_{N;s_0}$.\footnote{The partial ordering on terms comes from the monoid of bending parameters.  
See \cite[Section~3]{GHKK}.}
Since $\vartheta_p$ is an eigenfunction of the $H^{\times 3}$ action on $\ctal$, 
the other summands must have the same weight as $z^{\gv_{s_0} \lrp{\vartheta_p}}$ under this action.
Represent 
$\gv_{s_0} \lrp{\vartheta_p}$
pictorially in the following way, illustrated for $n=4$:

\input{g-vector.tex}

Note that
\eqn{
z^{\gv_{s_0} \lrp{\vartheta_p}} = \prod_{a+b+c = n} A_{\lrp{a,b,c}}^{g_{\lrp{a,b,c}}}.
}
Let $h_i = \diag\lrp{h_{i_1}, \dots, h_{i_n}} \in H$.
Then
\eqn{
\lrp{h_1,h_2,h_3} \cdot A_{\lrp{a,b,c}} = h_{1_1}\cdots h_{1_a} h_{2_1} \cdots h_{2_b} h_{3_1}\cdots h_{3_c} A_{\lrp{a,b,c}}.
}
Decompose $\lambda \in \chi^*\lrp{H}$ by $\lambda\lrp{h} = {h_1}^{\lambda_1} \cdots {h_n}^{\lambda_n}$.
Then the following picture lets us read off the $H^{\times 3}$ weight $\lrp{\abc}$ of 
$z^{\gv_{s_0} \lrp{\vartheta_p}}$, 
and in turn $\vartheta_p$ (denoted $w\lrp{\vartheta_p}$).

\input{Polytope3.tex}

There are two immediate consequences.
First,
{\prop{\label{prop:0sum}
Let $p \in \Xil\lrp{\Z^T}$ and write 
\eqn{
\gv_{s_0} \lrp{\vartheta_p} = \sum_{a+b+c = n} g_{\lrp{a,b,c}} e_{\lrp{a,b,c}}^*. 
}
Then $\vartheta_p$ is $\GL_n$-invariant (and hence $p\in \Xi \lrp{\Z^T}$),
if and only if
\eqn{
\sum_{a+b+c = n} g_{\lrp{a,b,c}} = 0. 
}
}}

Next,
{\prop{\label{prop:Pabc}
For each $\lrp{\abc} \in \chi^*\lrp{H^{\times 3}}$, define $P_{\abc} \subset \Xi$ to be the subset cut out by the hyperplanes described below:

\input{Polytope.tex}

Then $P_{\abc} \lrp{\Z^T}$ parametrizes a basis for $\lrp{\Vabc}^G$,
and $c^\gamma_{\alpha \beta} = \lrm{P_{\alpha,\beta,-w_0\lrp{\gamma}} \lrp{\Z^T}}$.
}}

We now have a representation theoretic interpretation of $\Xi$.
We'll pictorially represent the inequalities cutting out $\Xi$ in the initial seed as follows, bearing in mind that the sum of all entries must be 0.

\input{XiIneqs.tex}

Before finding the map $p^*:N\to M$, let's recall briefly the properties it must satisfy: \cite[p.~146]{GHK_birational} \footnote{Note that $M= M^\circ$ here-- all multipliers $d_i$ are $1$.}
\begin{enumerate}
        \item \label{p1} $\left.p^*\right|_{N_{\mathrm{uf}}} : n \mapsto \lrc{n, \cdot }$ and
        \item \label{p2} if $\pi: M \to M/N_{\mathrm{uf}}^\perp$ is the canonical projection, then $\pi\circ p^*:n \mapsto \lrb{\lrc{n, \cdot} : N_{\mathrm{uf}} \to \Z}$. 
\end{enumerate}

Here $N_{\mathrm{uf}}$ is the unfrozen sublattice of $N$-- the span of $\lrc{e_i}_{i \in I_{\mathrm{uf}}}$.  
Up to some map $N/N_{\uf} \to N_{\uf}^\perp $, $p^*$ is essentially just the exchange matrix $\epsilon$.
The ambiguity we can exploit-- the map $N/N_{\uf} \to N_{\uf}^\perp $-- amounts to choosing pairings between frozen indices, without considering any skew-symmetry requirements for this portion of our exchange matrix.

So we know what properties $p^*$ must satisfy, and we can compare Figures~\ref{fig:HiveBorder} and \ref{fig:Polytope} to further guide our efforts to write down a candidate $p^*$ map.
%
First, take $e_{\lrp{a,b,c}} \in N_{\mathrm{uf}}$. Using (\ref{p1}) we can immediately write down $p^* \lrp{e_{\lrp{a,b,c}}} = \lrc{e_{\lrp{a,b,c}}, \cdot }$.
For example, if we take  $e_{\lrp{2,1,1}}$:

\input{RHvectInt.tex}
then $p^* \lrp{e_{\lrp{2,1,1}}}$ must be

\input{pstarInt.tex}

Note that the sum of the entries of $p^* \lrp{e_{\lrp{2,1,1}}}$ is $0$, as are all sums indicated in Figure~\ref{fig:Polytope}.
So things look good so far.

Next, take $e_{\lrp{a,b,c}} \in N_{\mathrm{f}}$-- the frozen sublattice.
Decompose $p^*\lrp{e_{\lrp{a,b,c}} }$ as $p^*\lrp{e_{\lrp{a,b,c}} }_{\mathrm{uf}} +  p^*\lrp{e_{\lrp{a,b,c}} }_{\mathrm{f}} $.
Then (\ref{p2}) gives us the unfrozen portion $ p^*\lrp{e_{\lrp{a,b,c}} }_{\mathrm{uf}} $, and comparing Figures~\ref{fig:HiveBorder} and \ref{fig:Polytope} will suggest a candidate for the frozen portion.
Take for instance $e_{\lrp{3,1,0}}$:

\input{RHvector.tex}

Then (\ref{p2}) gives the unfrozen portion of $p^*\lrp{ e_{\lrp{3,1,0}} }$ as

\input{pstarvectUF.tex}

Now we use Figure~\ref{fig:HiveBorder} to find that $e_{\lrp{3,1,0}}$ has $\alpha = \lrp{1,-1,0,0}$, $\beta = \lrp{0,0,0,0}$, and $\gamma = \lrp{0,0,0,0}$.
So $p^*\lrp{ e_{\lrp{3,1,0}}}$ should have top entry $-\alpha_1 = -1$, entries of the top two rows summing to $-\alpha_2 = 1$, and all other sums from Figure~\ref{fig:Polytope} equal to $0$.
The obvious candidate then is

\input{pstarvect.tex}

Note again that the sum of the entries is $0$.

We can do the same procedure for every $e_{\lrp{a,b,c}}$.
The resulting map is given by

\input{TopEntry.tex}

\input{EdgeEntry.tex}

\input{InteriorEntry.tex}
and rotations of these.  Every entry not explicitly given is $0$.

So this is our proposed map $p^*$.  It certainly satisfies (\ref{p1}) and (\ref{p2}).
What we want to show now is that $p^*$ gives a unimodular equivalence between version (\ref{KTquot}) of the Knutson-Tao hive cone and $\Xi$ in the initial seed, 
and furthermore that $p^*W = \mathcal{W}_{GS}$-- so the representation theoretic Goncharov-Shen potential has a purely geometric description.

{\prop{\label{prop:WGS}\eqn{p^*W = \mathcal{W}_{GS}}}}
\begin{proof}
Consider the $\vartheta$-function $\vartheta_{\lrp{a,b,0}} = \lrp{\sum_{i= 0}^{n-a-1} z^{-\sum_{j=0}^i e_{\lrp{a,b-j,j}}}}$.
\eqn{
p^* \vartheta_{\lrp{a,b,0}}&= p^*\lrp{\sum_{i= 0}^{n-a-1} z^{-\sum_{j=0}^i e_{\lrp{a,b-j,j}}}}\\
&= \sum_{i= 0}^{n-a-1} z^{-\sum_{j=0}^i p^*\lrp{e_{\lrp{a,b-j,j}}}}\\
&= \sum_{i= 0}^{n-a-1} z^{-e_{\lrp{a,b-i,i}}^* + e_{\lrp{a-1,b-i,i+1}}^* -e_{\lrp{a,b-i-1,i+1}}^* + e_{\lrp{a+1,b-i-1,i}}^* }\\
&= \sum_{i=0}^{n-a-1} \frac{A_{\lrp{a-1,b-i,i+1}} A_{\lrp{a+1,b-i-1,i}}}{ A_{\lrp{a,b-i,i}} A_{\lrp{a,b-i-1,i+1}}}.
}
The last line above is expressed in \cite{GShen} as
\eqn{
\sum_{i= 0}^{n-a-1} \frac{\Delta_{a-1,b-i,i+1} \Delta_{a+1,b-i-1,i}}{\Delta_{a,b-i,i} \Delta_{a,b-i-1,i+1}}.
}
Summing over all $\vartheta$-functions in $W$ yields the potential $\mathcal{W}$ in \cite[Section~3.1]{GShen},
with each monomial summand corresponding to a different rhombus inequality.
\end{proof}

Pictorially, $p^*$ identifies the inequalities defining $\Xi$ in the initial seed with those defining the Knutson-Tao hive cone in the following way:

\input{KTcorrespondence.tex}

{\prop{\label{prop:unimodular}The map $p^*$ is unimodular, so in the initial seed $\Xi$ is unimodularly equivalent to the Knutson-Tao hive cone.}}

{\remark{Keep in mind that upon identifying our tropical spaces with real vector spaces, the domain of $p^*$ will look like $\R^{\mathcal{H}}/ \R \cdot {\bf{1}}_{\mathcal{H}}$ and the codomain will look like the subspace $V$ of $\R^{\mathcal{H}}$ in which the sum of all entries is 0.
Note that if we view the two copies of $\R^{\mathcal{H}}$ as dual spaces in the obvious way, then ${{\bf{1}}_{\mathcal{H}}}^\perp$ is exactly $V$, and so the domain and codomain of $p^*$ are also dual spaces.}}

\begin{proof}
First note that
\eq{
\sum_{\lrp{a,b,c} \in \mathcal{H}} p^* \lrp{e_{\lrp{a,b,c}}} = p^*\lrp{{\bf{1}}_{\mathcal{H}}} = 0.
}{eqn:sum0}

Next, I claim that
\eqn{
\Sp_\Z \lrc{e_{\lrp{n,0,0}}^*, p^*\lrp{e_{\lrp{a,b,c}}} }_{\lrp{a,b,c} \in \mathcal{H}} = \Z^{\mathcal{H}}.
}
On account of (\ref{eqn:sum0}), an immediate corollary of this claim would be that
\eqn{
\lrc{e_{\lrp{n,0,0}}^*, p^*\lrp{e_{\lrp{a,b,c}}} }_{{\text{All except one }}\lrp{a,b,c} \in \mathcal{H}}
}
is a basis for $\Z^{\mathcal{H}}$,
and so
\eqn{
\lrc{p^*\lrp{e_{\lrp{a,b,c}}} }_{{\text{All except one }}\lrp{a,b,c} \in \mathcal{H}}
}
would have to be a basis for $\Z^{\mathcal{H}} \bigcap V$.
Since
\eqn{
\lrc{e_{\lrp{a,b,c}} }_{{\text{All except one }}\lrp{a,b,c} \in \mathcal{H}}
}
is a basis for $\R^{\mathcal{H}}/ \R \cdot {\bf{1}}_{\mathcal{H}}$, this would establish unimodularity of $p^*$.
On to the claim.

As seen in the proof of Proposition~\ref{prop:WGS}, for each rhombus defining an inequality of the Knutson-Tao hive cone, we get a vector in the image of $p^*$ having $1$'s as the entries of the obtuse vertices, $-1$'s as the entries of the accute vertices, and $0$'s elsewhere.
In addition, $p^*\lrp{e_{\lrp{n,0,0}}} = e_{\lrp{n,0,0}}^* - e_{\lrp{n-1,0,1}}^* $, displayed in Figure~\ref{fig:pstarTop}.
Adding this to the vector we've associated to the top vertical rhombus just shifts its non-zero entries 1 position southeast, giving $e_{\lrp{n-1,1,0}}^* - e_{\lrp{n-2,1,1}}^* $.

\input{SoutheastShift.tex}

We can use the other vertical rhombi along the northeast border to shift these entries along the rest of the border, ending with the vector $e_{\lrp{1,n-1,0}}^* - e_{\lrp{0,n-1,1}}^*$.
Now take the image of the southeast corner: $p^*\lrp{e_{\lrp{0,n,0}}} = e_{\lrp{0,n,0}}^* - e_{\lrp{1,n-1,0}}^*$.
Adding this to our previous result of  $e_{\lrp{1,n-1,0}}^* - e_{\lrp{0,n-1,1}}^*$ gives $e_{\lrp{0,n,0}}^* - e_{\lrp{0,n-1,1}}^*$.

\input{AddSoutheastCorner.tex}

We can use the other collection of rhombi along the northeast border to shift the non-zero entries of this vector northwest along the border, starting with the rhombus containing the southeast corner $\lrp{0,n,0}$.

\input{NorthwestShift.tex}

So we've found two vectors in the image of $p^*$ to associate to each vertex $v_{\lrp{a,b,1}}$ along the $c=1$ line of $\mathcal{H}$:
$e_{\lrp{a+1,b,0}}^* - e_{\lrp{a,b,1}}^*$ oriented diagonally and
$e_{\lrp{a,b+1,0}}^* - e_{\lrp{a,b,1}}^*$ oriented horizontally.
Now introduce $e_{\lrp{n,0,0}}^*$.
Since $e_{\lrp{n,0,0}}^* - e_{\lrp{n-1,0,1}}^*$ and $e_{\lrp{n,0,0}}^*$ are in
$$\Lambda:= \Sp_\Z \lrc{e_{\lrp{n,0,0}}^*, p^*\lrp{e_{\lrp{a,b,c}}} }_{\lrp{a,b,c} \in \mathcal{H}},$$
so is $e_{\lrp{n-1,0,1}}^*$.
Next we use the other vector associated to $v_{\lrp{n-1,0,1}}$ (the horizontally oriented one $e_{\lrp{n-1,1,0}}^* - e_{\lrp{n-1,0,1}}^*$ this time) to see that
$e_{\lrp{n-1,1,0}}^*$ is also in $\Lambda$.
We then go one step southeast to the vertex $v_{\lrp{n-2,1,1}}$ and repeat this process, starting with the diagonally oriented vector and following up with the horizontally oriented vector,
to find that $e_{\lrp{n-2,1,1}}^*$ and $e_{\lrp{n-2,2,0}}^*$
are in $\Lambda$
as well.
Continuing southeast gives every $e_{\lrp{a,b,c}}^*$ with $c=0$ or $1$.
Now we'll push toward the southwest corner using rhombi of the only remaining orientation.
For each vertex $v_{\lrp{a,b,2}}$ along the $c=2$ line, there is a single rhombus having this as one of its vertices, two vertices with $c=1$, and one vertex with $c=0$.
Combining the vector associated to this rhombus with $\lrc{e_{\lrc{a,b,c}}^*}_{c=0 \text{ or } 1}$, we find that
$e_{\lrp{a,b,2}}^*$ is also in $\Lambda$.
Repeat this for $c=3$, then $4$, and so on out to $n$.
So each
$e_{\lrp{a,b,c}}^*$ is in $\Lambda$, and $\Lambda = \Z^{\mathcal{H}}$, completing the proof.
\end{proof}

\subsection{\label{subsec:Tk}Discussion of \texorpdfstring{$H^{\times 3}$}{H} action and the weight map}

Let $p^*_2$ be the composition $N \stackrel{p^*}{\longrightarrow} M \stackrel{\pi}{\longrightarrow} M/{N_{\uf}}^\perp$,
and denote the kernel of $p^*_2$ by $K$.
The inclusion $K \subset N$ induces an inclusion of tori $T_K \subset T_N$, and so an action of $T_K$ on $T_N$.
Furthermore, it induces a map $$T_M = \Spec\lrp{\kk \lrb{N}} \to T_K^* = \Spec \lrp{\kk\lrb{K}}.$$
Since $p^*$ commutes with mutation,
\begin{enumerate}
	\item it defines a map $p:\ctag \to \ctag^\vee$, 
	\item the action of $T_K$ on $T_N$ extends to an action on $\ctag = \bigcup_s T_{N;s}$, and
	\item \label{weight}it gives a map $\ctag^\vee = \bigcup_s T_{M;s} \to T_K^*$.
\end{enumerate}
This is discussed in greater detail in \cite[Section~2]{GHK_birational}.

Here I'll identify the action of $T_K$ on $\ctag$ with the $H^{\times 3}$ action, 
and the tropicalization of (\ref{weight}) with the map $w$ sending $p \in \ctag^\vee \lrp{\Z^T}$ to the $H^{\times 3}$-weight of the corresponding $\tf$-function $w\lrp{\vartheta_p}$.

The $H^{\times 3}$ action scales the decorations $\lrp{x_\bullet,y_\bullet,z_\bullet}$.
We decompose $h \in H^{\times 3}$ as $$h=\lrp{\lrp{h_{x_1},\dots, h_{x_n}},\lrp{h_{y_1},\dots, h_{y_n}},\lrp{h_{z_1},\dots, h_{z_n}}},$$
where, {\it{e.g.}}, $h_{x_i}$ is the scale factor for $x_i$.
Each component defines a one-parameter subgroup of $T_N$, which we'll show is in fact contained in $T_K$.
For instance, take $n=5$. 
Then the scaling coming from $h_{x_3}$ can be represented by

\input{hx3.tex}

For arbitrary $n$, 
$h_{x_i}$ corresponds to the cocharacter 
\eqn{
n_{x_i}:= \sum_{\lrp{a,b,c} \in \mathcal{H}, a \geq i} e_{\lrp{a,b,c}},
}
$h_{y_i}$ to 
\eqn{
n_{y_i}:= \sum_{\lrp{a,b,c} \in \mathcal{H}, b \geq i} e_{\lrp{a,b,c}},
}
and $h_{z_i}$ to 
\eqn{
n_{z_i}:= \sum_{\lrp{a,b,c} \in \mathcal{H}, c \geq i} e_{\lrp{a,b,c}}.
}

{\prop{
\eqn{
K = \Sp_{\Z} \lrc{ n_{x_i}, n_{y_i}, n_{z_i}}_{1\leq i \leq n}
}
In other words, the $H^{\times 3}$ action on $\ctag$ is precisely the $T_K$ action induced by the inclusion $K \subset N$. 
}}

\begin{proof}
We'll show that 
\eqn{
\Sp_{\Z} \lrc{p^*\lrp{ n_{x_i}}, p^*\lrp{ n_{y_i}}, p^*\lrp{ n_{z_i}}}_{1\leq i \leq n}  = N_\uf^\perp.
}
This implies containment in $K$, and unimodularity of $p^*$ boosts this containment to an equality.

We simply compute.
\eqn{
p^*\lrp{n_{x_i}}=& p^*\lrp{\sum_{\lrp{a,b,c} \in \mathcal{H}, a \geq i} e_{\lrp{a,b,c}}}\\
		=& p^*\lrp{e_{\lrp{n,0,0}}} + \tc{purple}{\sum_{a=i}^{n-1} p^*\lrp{e_{\lrp{a,0,n-a}}}} + \tc{blue}{\sum_{a=i}^{n-1} p^*\lrp{ \sum_{b=1}^{n-a} e_{\lrp{a,b,n-a-b}}}}\\
		=& e_{\lrp{n,0,0}}^* - e_{\lrp{n-1,0,1}}^* +\tc{purple}{ \sum_{a=i}^{n-1} e_{\lrp{a,0,n-a}}^* - e_{\lrp{a,1,n-a-1}}^* - e_{\lrp{a-1,0,n-a+1}}^*+ e_{\lrp{a-1,1,n-a}}^*  }\\
		 &+ \tc{blue}{\sum_{a=i}^{n-1} -e_{\lrp{a+1,0,n-a-1}}^*+ e_{\lrp{a,0,n-a}}^* + e_{\lrp{a,1,n-a-1}}^* -e_{\lrp{a-1,1,n-a}}^* }\\
		=& e_{\lrp{n,0,0}}^* - e_{\lrp{n-1,0,1}}^*+\tc{purple}{ e_{\lrp{n-1,0,1}}^* - e_{\lrp{n-1,1,0}}^* - e_{\lrp{i-1,0,n-i+1}}^*+ e_{\lrp{i-1,1,n-i}}^*  }\\
		 & \tc{blue}{ -e_{\lrp{n,0,0}}^*+ e_{\lrp{i,0,n-i}}^* + e_{\lrp{n-1,1,0}}^* -e_{\lrp{i-1,1,n-i}}^* }\\
		=& e_{\lrp{i,0,n-i}}^* - e_{\lrp{i-1,0,n-i+1}}^*
}
Similarly, 
\eqn{
p^*\lrp{n_{y_i}} = e_{\lrp{n-i, i, 0}}^* -e_{\lrp{n-i+1,i-1,0}}^*
}
and
\eqn{
p^*\lrp{n_{z_i}} = e_{\lrp{0,n-i, i}}^* -e_{\lrp{0,n-i+1,i-1}}^*.
}
Then clearly
\eqn{
\Sp_{\Z} \lrc{p^*\lrp{ n_{x_i}}, p^*\lrp{ n_{y_i}}, p^*\lrp{ n_{z_i}}}_{1\leq i \leq n}  = N_\uf^\perp,
}
(recall that $N = \Z^{\mathcal{H}}/{\bf{1}}_{\mathcal{H}} \cdot \Z $)
and so 
\eqn{
K = \Sp_{\Z} \lrc{ n_{x_i}, n_{y_i}, n_{z_i}}_{1\leq i \leq n}.
}
\end{proof}

To see that (\ref{weight}) tropicalizes to the weight map, first restrict to tori for a fixed seed $s$.
Then $\vartheta_p$ is a finite sum of characters on $T_{N;s}$, and since $\vartheta_p$ is an $H^{\times 3}$ eigenfunction, each of these characters has the same $H^{\times 3} = T_K$ weight.
Let one of the characters be $z^m$.  Then the weight of $\vartheta_p$ under the $T_K$ action is
the map $z^k \mapsto z^{\lra{k,m}}$ for $z^k \in T_K$.
In other words, the $T_K$ weight of $z^m$ is $\lrp{m \mod K^\perp} \in K^*$.\footnote{As a sanity check, note that in the scattering diagram approach of \cite{GHKK}, all scattering functions have the form $\lrp{1+z^{p^*(n)}}^c$ for some $n\in N_{\uf}$ and $c\in \Z_{>0}$.
As a result, the exponents of $\vartheta_p$'s summands differ by elements $p^*\lrp{N_\uf}$.
It is straightforward to check that $p^*\lrp{N_\uf} \subset K^\perp$, so all of the summands of $\vartheta_p$ do indeed have the same weight.
} 
The map ${m \mapsto m \mod K^\perp}$ dualizes the inclusion $K \hookrightarrow N$, so for each seed $s$ the weight of $\vartheta_p$ is the tropicalization of (\ref{weight}).
Since this holds when we restrict to every torus, it holds for all of $\ctag^\vee$.

As alluded to previously, there is a related action on $\cAp$.
Let $\widetilde{K}$ be the kernel of 
\eqn{
N\oplus M  &\to M/N_\uf^\perp\\
\lrp{n,m} &\mapsto p_2^*(n) - m.
}
The surjection $\pi:\cAp \to T_M$ is $T_{\widetilde{K}}$-equivariant.
The fact that the exchange matrix is full rank implies that $\pi$ is isomorphic to the trivial bundle $\ctag \times T_M$.\cite[Lemma~B.7]{GHKK}
This is used in \cite[Proof of Corollaries~0.20 and 0.21, page~581]{GHKK} to translate basis results for $\overline{\cA}_\prin$ to $\cta$. 

\subsection{Ray representation of \texorpdfstring{$\Xi$}{Xi}}

In \cite{Magee}, I showed that for the decorated flag variety, the cone $\Xi_{\Fld}$ is generated by the $\gv$-vectors for Pl\"ucker coordinates.
This is probably the most natural generating set for the homogeneous coordinate ring of the flag variety,
and I wonder if the generators of $\Xi$ can fill this role for $\ctb$.
This is necessarily vague and subjective.
I am primarily asking if the $\vartheta$-functions corresponding to generators of $\Xi$ in the initial seed have a simple explicit description.
In this subsection, I describe the rays generating $\Xi$ and give partial results relating these rays to functions.
Since $\Xi$ is not strictly convex for $G=\GL_n$, we'll temporarily restrict to $\SL_n$.

{\prop{\label{prop:frozen_edges} The $\gv$-vectors of all frozen variables generate edges of $\Xi$.}}

\begin{proof}
These $\gv$-vectors have a single non-zero entry.
The frozen variable $A_{\lrp{a,b,c}}$ has a 1 in the $\lrp{a,b,c}$ position, which is on a boundary edge of the triangle in Figure {\ref{fig:4triangulation}}.
For example, in $\SL_4$ the entries of the $\gv$-vector for $A_{\lrp{2,0,2}}$ are \\

\input{frozen_g.tex}

We can include the boxes and arrows representing $\Xi$'s defining inequalities: \\

\input{frozen_g_ineq.tex}

The line spanned by such a $\gv$-vector is the intersection of the hyperplanes
defined by the boxes for the remaining frozen variables and the arrows parallel to the boundary edge containing $v_{\lrp{a,b,c}}$.
So for $A_{\lrp{2,0,2}}$, we get the line described by the following picture:\\

\input{frozen_line.tex}

The ray $\R_{\geq 0} \cdot \gv\lrp{A_{\lrp{a,b,c}}}$ contained in this line satisfies the remaining inequalities and is an edge of $\Xi$.
\end{proof}

{\prop{\label{prop:s0edges} The $\gv$-vectors of all initial seed variables generate edges of $\Xi$.}}

\begin{proof}
Proposition {\ref{prop:frozen_edges}} took care of the frozen variables.
For the unfrozen variable $A_{\lrp{a,b,c}}$, there is a $1$ in the interior of our triangle with every other entry 0.
For starters, take the hyperplanes defined by each box.
Then fix a box and add as many consecutive arrows in the string emanating from it as possible without hitting $v_{\lrp{a,b,c}}$.
Doing this for all boxes gives enough hyperplanes to determine the line given by a free parameter in position $v_{\lrp{a,b,c}}$ and 0 elsewhere.
It really gives more hyperplanes than needed, but that's not a problem.
For example, for $A_{\lrp{2,1,1}}$ we would have the following picture:

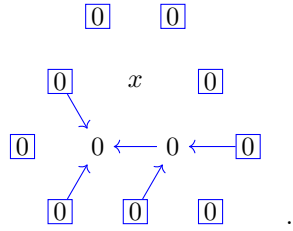
\input{s0_line.tex}

The ray $\R_{\geq 0} \cdot \gv\lrp{A_{\lrp{a,b,c}}}$ contained in this line satisfies the remaining inequalities and is an edge of $\Xi$.
\end{proof}

The functions associated to the remaining edges take more effort to describe.
We'll figure out what these edges actually are before worrying what functions they might correspond to.
Consider for a moment the inequalities represented by arrows.
These are described in Figure \ref{fig:XiIneqs}.
Staring at this for a little while suggests the following picture for some of the remaining rays of $\Xi$.
Take an interior vertex $v_{\lrp{a,b,c}}$.
There are three subquivers $Q_L$ (from Proposition \ref{prop:op}) starting at a frozen vertex and ending at $v_{\lrp{a,b,c}}$.
Draw a line segment through the vertices of each of these subquivers.
Here's an example:\\

\input{3vertex.tex}

Now put a ``$-1$'' in the $\lrp{a,b,c}$ position, a single ``$1$'' along each of the three segments, and ``$0$'' elsewhere.
The idea is that each of the three arrows pointing to $v_{\lrp{a,b,c}}$ indicate that the sum of the entries along one of the three segments should be non-negative.
If we make any entry negative, we'll get something outside of the cone generated by vectors in Propositions \ref{prop:frozen_edges} and \ref{prop:s0edges}.
The construction described is the simplest way to achieve this without violating any inequalities.\\

\input{3vertex_filled.tex}

Let's call such a picture a {\it{trivalent vertex}}.

{\prop{\label{prop:3vertex_edges} The vectors associated to trivalent vertices generate edges of $\Xi$.}}

\begin{proof}
Consider all of the hyperplanes defining faces of $\Xi$.
It is easiest to say which to exclude from our intersection.
Essentially, we want to intersect all of the hyperplanes associated to inequalities that should reduce to equalities for the vector in question, and only these hyperplanes.
So in Figure~\ref{fig:3vertex_filled}, we would remove exactly the blue boxes and arrows corresponding to strict inequalities, and the intersection of the remaining hyperplanes
is a line containing the given vector.
To do this for an arbitrary vector associated to a trivalent vertex,
we start by going to the position of the entry $1$ along each line segment.
If this isn't on the boundary, then the box at the end of the line segment and all arrows leading to this entry apart from the last one give equalities.
But the arrow whose tip hits this entry gives a strict inequality, as do the arrows coming after it, until we get to the arrow whose tip hits the $-1$ entry.
So given the segment

\input{segment.tex}

we would {\emph{not}} include the following hyperplanes in the intersection:

\input{segment_filled.tex}

Do this for all three line segments and then intersect all remaining hyperplanes.
The result is the real span of the given vector, and its $\R_{\geq 0}$ span also satisfies the inequalities that have been omitted from the intersection.
\end{proof}

The next thing to notice is that we can overlay two trivalent vertices, and as long as none of the line segments are colinear we'll have two line segments intersecting at a vertex.
If we make the entry of this vertex $1$, we'll
get a vector outside the span of the vectors previously described.
For example, take

\input{3vertexOverlay.tex}

which gives the vectors

\input{Overlaid3vertexVector.tex}

and

\input{Overlaid3vertexVector2.tex}

Consider either of these two vectors.
If it were the sum of vectors described previously, we'd have to take at least two vectors associated to trivalent vertices to account for the two minus signs.
Then the sum of all entries must be at least $4$, but it is in fact only $3$.
So it is indeed outside of the span of the vectors described previously, and it clearly lies in our cone.
The proof that it generates an edge of $\Xi$ is basically identical to the proof of Proposition~\ref{prop:3vertex_edges}.

Next, there is no reason to limit ourselves to only overlaying two trivalent vertices.
We can overlay as many as we want.
Say we overlay $k$ of them.
Then we are describing vectors for which $k$ entries are $-1$.
As long as we place our $1$'s in such a way that our vector cannot be a positive combination the vectors associated to $k-1$ or fewer trivalent vertices,
we will get an edge of $\Xi$ by the argument used above.
In particular, if we ensure that the sum of the entries is less than could be achieved with such positive combinations, our vector must be an edge of $\Xi$.

{\prop{Every edge of $\Xi$ is generated by an initial seed $\gv$-vector, a trivalent vertex, or overlaid trivalent vertices.}}

\begin{proof}
The initial seed $\gv$-vectors already generate the entire positive orthant,
so any additional edge of $\Xi$ must have some negative entry.
Negative entries must be at interior vertices on account of the box inequalities.
Say the entry at $\lrp{a,b,c}$ is negative, with value $-x$.
Then the three incoming arrows at $v_{\lrp{a,b,c}}$ indicate that the sum of the remaining entries along each of the three line segments leading to $v_{\lrp{a,b,c}}$ must be at least $x$.
However, if the sum is more than $x$, we would be able to realize this vector as a sum of vectors from trivalent vertices and vectors in the positive orthant.
You can see this by restricting to a line segment first.
It's easy to see for this restriction, and the result transfers over directly as the trivalent vertices are made up of three line segments, 
and the position of a $1$ along one line segment is completely independent of the other two line segments.
It follows that if some edge of $\Xi$ lies outside of the positive orthant, it must be generated by a trivalent vertex or a collection of overlaid trivalent vertices.
So this is indeed all of the edges of $\Xi$.
\end{proof}

{\prop{\label{prop:gMGS}Consider the trivalent vertex $p$ having $1$'s in positions $\lrp{a_1,b_1,c_1}$, $\lrp{a_2,b_2,c_2}$, and $\lrp{a_3,b_3,c_3}$, 
labeled such that the top left $1$ is in position $\lrp{a_1,b_1,c_1}$, the rightmost $1$ is at $\lrp{a_2,b_2,c_2}$, and the bottom left $1$ is $\lrp{a_3,b_3,c_3}$.
Now take the triangle $\triangle$ with sides $b=b_1+1$, $c=c_2+1$, and $a=a_3+1$.
Then $\vartheta_p$ is obtained by performing the maximal green sequence of Proposition~\ref{prop:MGS} on the subquiver with vertices $\triangle$.
}}

Let's illustrate the claim first.  
Suppose we take 

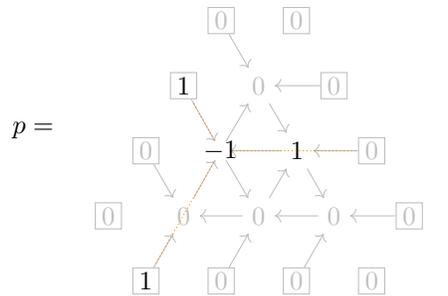
\input{g-tri.tex}

Then $\vartheta_p$ is produced by the mutation sequence

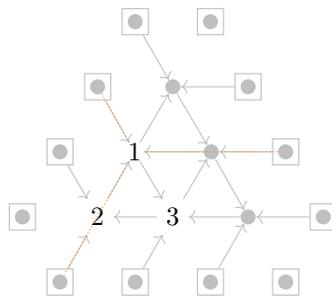
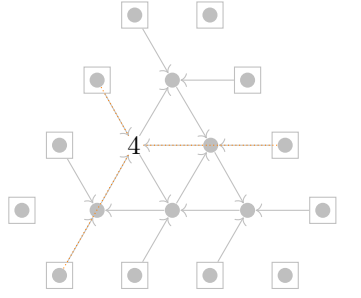
\input{g-tri_mut.tex}
then
\input{g-tri_mut2.tex}

\begin{proof}
Each vertex of mutation in this sequence is green, so when mutating at $v_k$ the terms coming from arrows emanating from $v_k$ vanish on the central fiber of $\cA_{\prin,s_0}$.
Only those arrows terminating at $v_k$ contribute to the $\gv$-vector.
Let 
$A_{\lrp{a,b,c}_k}$ 
be the cluster variable obtained after the $k^{\mathrm{th}}$ mutation at $v_{\lrp{a,b,c}}$ in the sequence.
Then, using the quivers from the proof of Proposition~\ref{prop:MGS}, 
$\gv_{s_0} \lrp{A_{\lrp{a,b,c}_k}} $
is given by 
\eqn{
	\gv_{s_0} \lrp{ A_{\lrp{n -b_1 -c_2-k, b_1, c_2+k}}} + \gv_{s_0} \lrp{A_{\lrp{a-k+1,b+k,c_2}}} + \gv_{s_0} \lrp{A_{\lrp{a-1,b,c+1}_{k-1}}} -\gv_{s_0} \lrp{A_{\lrp{a,b,c}_{k-1}}} 
}
if $ c = c_2+1$,
and
\eqn{
	\gv_{s_0} \lrp{A_{\lrp{a+1,b,c-1}_k}} + \gv_{s_0} \lrp{A_{\lrp{a-1,b,c+1}_{k-1}}} -\gv_{s_0} \lrp{A_{\lrp{a,b,c}_{k-1}}} 
}
otherwise.
I claim that the for each $v_{\lrp{a,b,c}} \in \triangle$, 
\eqn{
\gv_{s_0} \lrp{A_{\lrp{a,b,c}_k}} = e_{\lrp{n-b_1-c_2-k, b_1, c_2+k }}^* + e_{\lrp{n-b-c_2-k,b+k,c_2}}^* + e_{\lrp{a-k,b,c+k}}^* - e_{\lrp{n-b-c_2-k,b,c_2+k}}^* .
}
For $k=1$, 
$\gv_{s_0} \lrp{A_{\lrp{a,b,c}_1}} $
is given by 
\eqn{
	e_{\lrp{n -b_1 -c_2-1, b_1, c_2+1 }}^* + e_{\lrp{a,b+1, c_2}}^* + e_{\lrp{a-1,b,c+1}}^* - e_{\lrp{a,b,c}}^* 
}
if $ c = c_2+1$, in agreement with the claim,
and
\eqn{
	\gv_{s_0} \lrp{A_{\lrp{a+1,b,c-1}_1}} + e_{\lrp{a-1,b,c+1}}^* - e_{\lrp{a,b,c}}^* 
}
otherwise.
When $c-1 = c_2+1$, we have 
\eqn{
\gv_{s_0} \lrp{A_{\lrp{a+1,b,c_2+1}_1}} =
	e_{\lrp{n -b_1 -c_2-1, b_1, c_2+1 }}^* + e_{\lrp{a+1,b+1, c_2}}^* + e_{\lrp{a,b,c}}^* - e_{\lrp{a+1,b,c_2+1}}^*,
}
so
\eqn{
	\gv_{s_0} \lrp{A_{\lrp{a,b,c_2+2}_1}}= e_{\lrp{n -b_1 -c_2-1, b_1, c_2+1 }}^* + e_{\lrp{a+1,b+1, c_2}}^*  + e_{\lrp{a-1,b,c_2+3}}^* - e_{\lrp{a+1,b,c_2+1}}^*,
}
which again agrees with the claim.
Assume the claim holds for all $c' < c$.
Then 
\eqn{
	\gv_{s_0} \lrp{A_{\lrp{a,b,c}_1}} =& \lrp{  e_{\lrp{n-b_1-c_2-1, b_1, c_2+1 }}^* + e_{\lrp{n-b-c_2-1,b+1,c_2}}^* + e_{\lrp{a,b,c}}^* - e_{\lrp{n-b-c_2-1,b,c_2+1}}^*  }\\
    					   &+ e_{\lrp{a-1,b,c+1}}^* - e_{\lrp{a,b,c}}^* \\
					  =&  e_{\lrp{n-b_1-c_2-1, b_1, c_2+1 }}^* + e_{\lrp{n-b-c_2-1,b+1,c_2}}^*  + e_{\lrp{a-1,b,c+1}}^* - e_{\lrp{n-b-c_2-1,b,c_2+1}}^*, 
}	
which proves the claim for $k=1$.
Now suppose it holds for all $k' <k$.
Then 
$\gv_{s_0} \lrp{A_{\lrp{a,b,c_2+1}_k}} $ is given by
\eqn{
	&e_{\lrp{n -b_1 -c_2-k, b_1, c_2+k}}^* + e_{\lrp{a-k+1,b+k,c_2}}^* \\
	&+ \lrp{  e_{\lrp{n-b_1-c_2-k+1, b_1, c_2+k-1 }}^* + e_{\lrp{n-b-c_2-k+1,b+k-1,c_2}}^* + e_{\lrp{a-k,b,c_2+1+k}}^* -     e_{\lrp{n-b-c_2-k+1,b,c_2+k-1}}^*  }\\
	&- \lrp{  e_{\lrp{n-b_1-c_2-k+1, b_1, c_2+k-1 }}^* + e_{\lrp{n-b-c_2-k+1,b+k-1,c_2}}^* + e_{\lrp{a-k+1,b,c_2+1+k-1}}^* - e_{\lrp{n-b-c_2-k+1,b,c_2+k-1}}^*  }\\
	=& e_{\lrp{n -b_1 -c_2-k, b_1, c_2+k}}^* + e_{\lrp{a-k+1,b+k,c_2}}^* + e_{\lrp{a-k,b,c_2+1+k}}^* - e_{\lrp{a-k+1,b,c_2+k}}^*,
}
in agreement with the claim,
and $\gv_{s_0} \lrp{A_{\lrp{a,b,c}_k}} $, $c \neq c_2+1$, is given by
\eqn{
	&\gv_{s_0} \lrp{A_{\lrp{a+1,b,c-1}_k}}\\
	&+ \lrp{  e_{\lrp{n-b_1-c_2-k+1, b_1, c_2+k-1 }}^* + e_{\lrp{n-b-c_2-k+1,b+k-1,c_2}}^* + e_{\lrp{a-k,b,c+k}}^* - e_{\lrp{n-b-c_2-k+1,b,c_2+k-1}}^*  }\\
	&- \lrp{  e_{\lrp{n-b_1-c_2-k+1, b_1, c_2+k-1 }}^* + e_{\lrp{n-b-c_2-k+1,b+k-1,c_2}}^* + e_{\lrp{a-k+1,b,c+k-1}}^* - e_{\lrp{n-b-c_2-k+1,b,c_2+k-1}}^*  }\\
	=& \gv_{s_0} \lrp{A_{\lrp{a+1,b,c-1}_k}} + e_{\lrp{a-k,b,c+k}}^* -e_{\lrp{a-k+1,b,c+k-1}}^* .
}
As before, when $c-1=c_2+1$, we have
\eqn{
\gv_{s_0} \lrp{A_{\lrp{a+1,b,c_2+1}_k}} =
	e_{\lrp{n -b_1 -c_2-k, b_1, c_2+k }}^* + e_{\lrp{a-k+2,b+k, c_2}}^* + e_{\lrp{a+1-k,b,c_2+1+k}}^* - e_{\lrp{a-k+2,b,c_2+k}}^*,
}
which agrees with the claim.
So suppose it holds for all $c'<c$. Then
\eqn{
	\gv_{s_0} \lrp{A_{\lrp{a,b,c}_k}} =& \gv_{s_0} \lrp{A_{\lrp{a+1,b,c-1}_k}} + \gv_{s_0} \lrp{A_{\lrp{a-1,b,c+1}_{k-1}}} -\gv_{s_0} \lrp{A_{\lrp{a,b,c}_{k-1}}} \\
					  =&  \lrp{ e_{\lrp{n-b_1-c_2-k, b_1, c_2+k }}^* + e_{\lrp{n-b-c_2-k,b+k,c_2}}^*         + e_{\lrp{a-k+1,b,c+k-1}}^* - e_{\lrp{n-b-c_2-k,b,c_2+k}}^* }\\
					   &+ \lrp{ e_{\lrp{n-b_1-c_2-k+1, b_1, c_2+k-1 }}^* + e_{\lrp{n-b-c_2-k+1,b+k-1,c_2}}^* + e_{\lrp{a-k,b,c+k}}^* -     e_{\lrp{n-b-c_2-k+1,b,c_2+k-1}}^* }\\
					   &- \lrp{ e_{\lrp{n-b_1-c_2-k+1, b_1, c_2+k-1 }}^* + e_{\lrp{n-b-c_2-k+1,b+k-1,c_2}}^* + e_{\lrp{a-k+1,b,c+k-1}}^* - e_{\lrp{n-b-c_2-k+1,b,c_2+k-1}}^* }\\
    					  =& e_{\lrp{n-b_1-c_2-k, b_1, c_2+k }}^* + e_{\lrp{n-b-c_2-k,b+k,c_2}}^* + e_{\lrp{a-k,b,c+k}}^* - e_{\lrp{n-b-c_2-k,b,c_2+k}}^* .
}	
This proves the claim.

Now take $k= c_1-c_2$,
and take $\lrp{a,b,c} = \lrp{a_3+c_1-c_2,b_3,c_2+a_2-a_3}$.
Then
\eqn{
	\gv_{s_0} \lrp{A_{\lrp{a,b,c}_k}} =& e_{\lrp{a_1, b_1, c_1}}^* + e_{\lrp{n-b_3-c_1, b_3+c_1-c_2,c_2}}^*+ e_{\lrp{a_3,b_3,a_2-a_3+c_1}}^* - e_{\lrp{n-b_3-c_1,b_3,c_1}}^* \\
					  =& e_{\lrp{a_1, b_1, c_1}}^* + e_{\lrp{a_2, b_2,c_2}}^* + e_{\lrp{a_3,b_3,c_3}}^* - e_{\lrp{a_2,b_3,c_1}}^*, 
}	
proving the proposition.
\end{proof}

Taking $\triangle =\trir{n-2}$ leads to a nice observation about the maximal green sequence of Proposition~\ref{prop:MGS}.
Call the final seed of the sequence $s$.

{\cor{\label{cor:MGSDualizes}Let $\mumgs$ be the automorphism of $\cta$ induced by the maximal green sequence of Proposition~\ref{prop:MGS}, along with the indicated permutation of frozen vertices.
Then the assignment 
$$A_{\lrp{a,b,c}_{s_0}} \mapsto \mumgs^*\lrp{A_{\lrp{a,b,c}_s}}$$ 
induces an automorphism of $\ssO\lrp{\cta}$ sending $\lrp{\Vabc}^G$ to $\lrp{V_\alpha^* \otimes V_\gamma^* \otimes V_\beta^*}^G$. 
}}

\begin{proof}
Since the only difference between $Q_s$ and $Q_{s_0}$ is the overall orientation, the variables of $s$ and $s_0$ have the same relations--
given the relation $$r\lrp{A_{\lrp{a_1,b_1,c_1}_{s_0}}, \dots, A_{\lrp{a_i,b_i,c_i}_{s_0}}} =1,$$ 
the relation
$$r\lrp{\mumgs^*\lrp{A_{\lrp{a_1,b_1,c_1}_{s}}}, \cdots,\mumgs^*\lrp{ A_{\lrp{a_i,b_i,c_i}_{s}} }} =1$$
must hold as well.
By Theorem~\ref{thm:FG}, 
every $\vartheta_p \in \Xi\lrp{\Z^T}$ is a Laurent polynomial in either of these two collections of variables.
Since $\Xi\lrp{\Z^T}$ generates $\ssO\lrp{\cta}$ and the relations between the two (identical but reordered) generating sets $\Xi \lrp{\Z^T}$ match,
this assignment gives an automorphism of $\ssO\lrp{\cta}$.

Now the claim is that if some $\vartheta_p$ has $H^{\times 3}$-weight $\lrp{\abc}$, 
then its image $\vartheta_p'$ under this automorphism
has $H^{\times 3}$-weight 
$\lrp{-w_0 \lrp{\alpha}, -w_0 \lrp{\gamma}, -w_0\lrp{\beta}}$.
It's sufficient to show that this holds for the cluster variables of $s_0$.
The weight of 
$A_{\lrp{a,b,c}_{s_0}}$
is 
\eqn{
\lrp{\lrp{h_{x_1},\dots,h_{x_{n-1}}}, \lrp{h_{y_1},\dots,h_{y_{n-1}}}, \lrp{h_{z_1},\dots,h_{z_{n-1}}}} \mapsto h_{x_1}\cdots h_{x_a} h_{y_1} \cdots h_{y_b} h_{z_1}\cdots h_{z_c}.
}
For the weight of 
$\mumgs^*\lrp{A_{\lrp{a,b,c}_s}}$, 
take $k=a$ and $b_1 = c_2 = a_3 = 0$ in the proof of Proposition~\ref{prop:gMGS}.
This yields
\eqn{
\gv_{s_0}\lrp{A_{\lrp{a,b,c}_s}} &= e_{\lrp{n-a,0,a}}^* + e_{\lrp{n-b-a,b+a,0}}^* + e_{\lrp{0,b,c+a}}^* - e_{\lrp{n-b-a,b,a}}^*\\
				 &= e_{\lrp{n-a,0,a}}^* + e_{\lrp{c,n-c,0}}^* + e_{\lrp{0,b,n-b}}^* - e_{\lrp{c,b,a}}^*,
}
so the weight is
\eqn{
\lrp{\lrp{h_{x_1},\dots,h_{x_{n-1}}}, \lrp{h_{y_1},\dots,h_{y_{n-1}}}, \lrp{h_{z_1},\dots,h_{z_{n-1}}}} \mapsto h_{x_1}\cdots h_{x_{n-a}} h_{y_1} \cdots h_{y_{n-c}} h_{z_1}\cdots h_{z_{n-b}}.
}
This proves the claim.\footnote{Recall that $G= \SL_n$ here.}
\end{proof}

It isn't clear yet whether the rays of $\Xi$ correspond to a ``simple'' collection of functions, but I hope that observations in this subsection provide a foundation for addressing this question.

\section{Recovering \texorpdfstring{$\Fld$}{decorated flag variety} and \texorpdfstring{$U$}{U} from \texorpdfstring{$\cta$}{configurations of flags}}

The decorated flag variety $\Fld$ is a partial minimal model for the double Bruhat cell $G^{e,w_0}$.
Fix $B_+, B_- \subset G$ to be the subgroups of upper and lower triangular matrices,
and take $V_\bullet^+$ and $V_\bullet^-$ to be their fixed flags.
Then $G^{e,w_0} \subset \Fld$ is the subset whose underlying flags $F_\bullet$ intersect both $V_\bullet^+$ and $V_\bullet^-$ generically.
That is, $F = \lrp{F_\bullet, f_\bullet}$ is in $G^{e,w_0}$ if and only if
each subspace $F_i$ intersects both $V_{n-i}^+$ and $V_{n-i}^-$ transversely.

This description, while satisfyingly simple, involves a choice-- fixing the pair $\lrp{B_+,B_-}$.
It would be philosophically more appealing to have a description that avoids such choices.
So, instead of choosing a single pair, we'll choose all pairs at once, and later we'll mod out to identify all of these choices.
First notice that $G=\GL\lrp{V}$ acts freely and transitively on the generic locus $\lrp{\Fld \times \Fl}^\times$ of $\Fld \times \Fl$.
It's clear that $G$ acts freely and transitively on ordered bases for $V$, and
the correspondence between ordered bases and generic pairs $\lrp{\lrp{X_\bullet,x_\bullet}, Y_\bullet}$ is straightforward.
Given an ordered basis $\lrp{v_1,\dots,v_n}$,
set $X_i:= \Sp\lrc{v_1,\dots,v_i}$, 
$Y_{n-i} := \Sp \lrc{v_{i+1},\dots, v_n}$,
and 
$x_i := v_i \mod X_{i-1}$.
On the other hand, given a generic pair $\lrp{\lrp{X_\bullet,x_\bullet}, Y_\bullet}$,
note that $x_{i+1}$ can be identified with an $i$-dimensional affine subspace of $X_{i+1}$: $x_{i+1} \text{``$=$''} \lrc{v \in X_{i+1} \left| v \mod X_i = x_{i+1} \right.}$.
Then $x_{i+1}$ and $Y_{n-i}$ intersect in a point since $X_i$ and $Y_{n-i}$ do, and $x_{i+1}$ is just a translation of $X_i$.
Set $v_{i+1} := x_{i+1} \bigcap Y_{n-i}$.
These two maps are clearly inverses of each other.
Now we can view $\Fld$ as the subset of $\Conf\lrp{\Fld,\Fld,\Fl}$ with $\lrp{A_1,B_3}$ generic, and $G^{e,w_0}$ as $\Conf\lrp{\Fld,\Fld,\Fl}^\times$.
The above description is based on \cite{GShen}.

So we start with $\ctag^\vee\lrp{\Z^T}$ and then take the slice whose $H_z$ weight $\gamma$ is $0$.\footnote{%
As in Remark~\ref{rem:CS}, using this condition we could recover the usual cluster structure for this space by defining a new collection of variables, say $\overline{A}_{\lrp{a,b,c}}:= A_{\lrp{a,b,c}}/A_{\lrp{0,n-c,c}}$.
But the point here is to avoid making choices, so we won't do that.}
This gives us a basis for $\ssO\lrp{G^{e,w_0}}$.
When we partially compactify to $\Fld$, we still ask for the first and third flags to intersect generically, so we are leaving off the divisors $D_{\lrp{i,0,n-i}}$.
Then in the initial seed, the Landau-Ginzburg potential in this case is 

\eqn{
W_{\Fld}= \sum_{a+b = n} \vartheta_{\lrp{a,b,0}} 
+ \sum_{b+c = n} \vartheta_{\lrp{0,b,c}}, 
}
where $ a,b,c \in \Z_{>0}$,  
\eqn{
\vartheta_{\lrp{a,b,0}}= 
\sum_{i= 0}^{n-a-1} z^{-\sum_{j=0}^i e_{\lrp{a,b-j,j}}}, 
}
and
\eqn{
\vartheta_{\lrp{0,b,c}}= 
\sum_{i= 0}^{n-b-1} z^{-\sum_{j=0}^i e_{\lrp{j,b,c-j}}}.
}
It is immediate from Proposition~\ref{prop:WGS} that $W_{\Fld}$ pulls back to the potential of \cite{GShen} on $\Conf\lrp{\Fld,\Fld,\Fl}$.
Furthermore, Proposition~\ref{prop:unimodular} and \cite[Theorem~3.2]{GShen} immediately imply that $\Xi_{\Fld}$ is unimodularly equivalent to the Gelfand-Tsetlin cone.
The polytope in $\Xi_{\Fld}$ where $\beta = -w_0 \lrp{\lambda}$ parametrizes a canonical basis for the irreducible representation $V_\lambda$.
{\remark{\label{rem:GT}
There is one semantic caveat worth mentioning here.  
I have generally seen the term ``Gelfand-Tsetlin cone'' applied to a strictly convex cone encoding polynomial $\GL_n$ representations.
$\ssO\lrp{\Fld}$ decomposes as a sum of rational $\GL_n$ representations, which include duals of polynomial representations.
Essentially, $\det$ is an invertible function on $\GL_n$ and this gives $\Xi_{\Fld}$ a 1-dimensional linear subspace.
That said, the cone defined by Gelfand and Tsetlin in \cite[Equation~3]{GT} encodes rational representations and is the cone identified with $\Xi_{\Fld}$ by $p^*$.
}}

Using \cite[Figure~31]{GShen}, the inequalities defining $\Xi_{\Fld}$ and the Gelfand-Tsetlin cone are identified via $p^*$ as follows:

\input{XiGT.tex}

Goncharov and Shen describe their potential as part of a 6-tuple defining a {\it{positive decorated geometric crystal}}.  
In this setting, $p^*W= \mathcal{W}_{\mathrm{GS}}$ plays the role of Berenstein and Kazhdan's potential $f$ from \cite{BKaz_lecture}.
See the appendix in the arXiv version of \cite{GShen} for details.

To describe $U$ in this way, we view $\Fl$ as the subset of $\Conf\lrp{\Fld,\Fl,\Fl}$ with $\lrp{A_1,B_3}$ generic, 
$U$ as the subset where $\lrp{A_1,B_2}$
is also generic,
and the cluster variety $\mathring{U}$ in $U$ as $\Conf^\times\lrp{\Fld,\Fl,\Fl}$.
Usually $\mathring{U}$ is described as the subset of $U$ (upper triangular unipotent matrices) where the minors $\Delta^{1,\dots,i}_{n-i+1,\dots,n}$ are non-vanishing.
Here we are getting to $U$ from $\Fl$ by requiring $\lrp{A_1,B_2}$ to be generic, 
and we get to $\mathring{U}$ from $U$ by requiring $\lrp{B_2,B_3}$ to be generic.

A basis for $\ssO\lrp{\mathring{U}}$ is given by taking the slice of $\ctag\lrp{\Z^T}$ with $H_y \times H_z$ weight $\lrp{\beta, \gamma} = 0$.
When we partially compactify to $U$, the divisors that we add are $D_{\lrp{0,i,n-i}}$.
The corresponding inequalities are the solid (as opposed to dashed) boxes and arrows of Figure~\ref{fig:XiGT}.
Then $\Xi_U$ is a simplicial cone of dimension $\binom{n}{2}$. 

%
%
%
\bibliography{bibliography}
\bibliographystyle{hep}

{\sc{School of Mathematics, University of Birmingham, Edgbaston, Birmingham B15 2TT, UK}}\\
{\it{e-mail:}} \href{mailto:t.magee@bham.ac.uk}{t.magee@bham.ac.uk}

\end{document}